\documentclass[reqno,11pt]{amsart}
\usepackage{epsfig}\usepackage{amsfonts,amsmath,amssymb,mathrsfs,verbatim,amsthm,amscd}
\numberwithin{equation}{section}
\newcommand{\beq}{\begin{equation}}
\newcommand{\ee}{\end{equation}}
\newcommand{\bea}{\begin{eqnarray}}
\newcommand{\eea}{\end{eqnarray}}

\usepackage{hyperref}
\def\stackreb#1#2{\ \mathrel{\mathop{#1}\limits_{#2}}}

\newcommand{\CC}{\mathbb C}
\newcommand{\RR}{\mathbb R}
\newcommand{\ZZ}{\mathbb Z}

\begin{document}

\vspace*{-2em}

\title[Infinite elliptic hypergeometric series]
{Infinite elliptic hypergeometric series:
\\
 convergence and difference equations%
}

\medskip

\author{D.\, I.~ Krotkov}
\address{National Research University Higher School of Economics, Moscow, Russia;
e-mail: dikrotkov@hse.ru}

\author{V.\, P.~ Spiridonov}%

\address{Laboratory of Theoretical Physics, JINR, Dubna, Moscow region, Russia and
National Research University Higher School of Economics, Moscow, Russia;
e-mail: spiridon@theor.jinr.ru
}

\maketitle

\vspace*{-1.0em}

\begin{abstract}
We derive finite difference equations of infinite order for theta hy\-per\-geo\-met\-ric series
and investigate the space of their solutions. In general, such infinite series
diverge, we describe some constraints on the parameters when they do converge. In particular,
we lift the Hardy-Littlewood criterion on the convergence of $q$-hy\-per\-geo\-met\-ric series
for $|q|=1, \, q^n\neq 1$, to the elliptic level and prove convergence of the infinite
${}_{r+1}V_r$ very-well poised elliptic hy\-per\-geo\-met\-ric series for restricted
values of $q$.
\end{abstract}


\tableofcontents

\section{Introduction}

Hypergeometric functions have very long history \cite{aar}. The notion of hy\-per\-geo\-met\-ric series
was introduced even before Newton proved the binomial theorem for infinite $_1F_0$-series.
The systematic study of such functions was launched by Euler, who introduced the
gamma function, computed the beta integral in terms of the gamma function, defined the key $_2F_1$-series
and considered the second order differential equation for that. He introduced also
the simplest $q$-hy\-per\-geo\-met\-ric series and found a number of exact identities
for them. A $q$-analogue of the hy\-per\-geo\-met\-ric equation was considered in
the middle of nineteenth century by Heine \cite{GR}. The theory of hy\-per\-geo\-met\-ric
functions was developing further on for a very long time in these two instances, which
have found very many applications in various fields of mathematics and theoretical physics.

At the turn of the Millenium a new type of hy\-per\-geo\-met\-ric functions
related to elliptic functions has been discovered, for a survey see \cite{spi:essays}.
It appeared that properties of elliptic functions explain some of the known beauties
of plain and $q$-hy\-per\-geo\-met\-ric functions. Moreover, at the elliptic level there emerged new
symmetries and many formulas became much more compact than their previously known degenerate
cases. By now, most of the classical constructions for the plain hy\-per\-geo\-met\-ric functions
have been uplifted to the elliptic level both at the univariate and multivariate settings.
Major part of the applications uses the elliptic hy\-per\-geo\-met\-ric series which satisfy by
definition a balancing condition imposed on the parameters \cite{spi:theta1}.
Some particular elliptic hy\-per\-geo\-met\-ric series emerge in the applications
to the Yang-Baxter equation \cite{FT}, in the theory of biorthogonal rational functions
\cite{spi:theta2,spi:essays} and
the theory of orthogonal polynomials with the spectrum dense on the unit circle \cite{Z1,Z2}.

However, some of the standard constructions are still waiting for proper un\-der\-stan\-ding at
the elliptic level. In the present paper we consider one of such problems ---
a straightforward generalization of the higher order hy\-per\-geo\-met\-ric differential equations
to finite difference equations of infinite order for theta hy\-per\-geo\-met\-ric series. Additionally, we
investigate convergence of infinite univariate theta hy\-per\-geo\-met\-ric and elliptic hy\-per\-geo\-met\-ric
series for particular restrictions on their parameters. The latter analysis resembles
the one performed by Hardy and Littlewood \cite{HL} for a particular
$q$-hy\-per\-geo\-met\-ric series, when $q=e^{2\pi i\chi}$
with a fixed irrational number $\chi$. A similar convergence problem for infinite series
arises at the elliptic level for arbitrary choice of parameters,
in sharp difference from the $q$-hy\-per\-geo\-met\-ric case.
We generalize the Hardy-Littlewood criterion to the infinite elliptic hy\-per\-geo\-met\-ric series
and prove convergence of such series in some particular cases. Especially, we show that for some choice
of parameters infinite well-poised elliptic hy\-per\-geo\-met\-ric series ${}_{r+1}V_r$ have the radius of
convergence bigger than 1.

\section{Plain hy\-per\-geo\-met\-ric series and their $q$-analogues}

Consider the standard $p$-Pochhammer symbol
\begin{equation}
(a; p)_{\infty} =\prod_{n=0}^\infty (1-ap^n), \qquad (a; p)_s
= \frac{(a;p)_\infty}{(ap^s; p)_\infty}, \quad |p|<1,
\end{equation}
and recall the following conventional notation for its $k$-fold product
$$
(a_1, \ldots ,a_k; p)_{\infty}  :=  (a_1; p)_\infty \cdot \ldots \cdot (a_k;p)_\infty.
$$
The Jacobi theta function we use is defined as follows
\begin{equation}
\theta(z; p)  :=  (z, pz^{-1} ; p)_\infty, \quad z\in\CC^\times,
\end{equation}
and the classical Jacobi triple-product identity gives the following expansion
\begin{equation}
\theta(z; p) = \frac{1}{(p;p)_\infty}\sum_{n \in \mathbb{Z}}p^{\binom{n}{2}}(-z)^n ,
\quad \binom{n}{2}=\frac{n(n-1)}{2}.
\label{TPI}\end{equation}
The key reflection and quasiperiodicity properties of this theta function are
\begin{equation}
\theta(z^{-1}; p) =\theta(pz; p) =-z^{-1}\theta(z; p).
\label{trafo}\end{equation}
As a consequence, we have
$$
\theta(p^kz; p) =(-z)^{-k}p^{-\binom{k}{2}}\theta(z; p), \quad k\in\ZZ.
$$

For an integer $n\in \mathbb{Z}_{\geq0}$, we define the elliptic shifted factorial,
or the elliptic Pochhammer symbol
\begin{equation}
\theta(t; p; q)_n  :=  \prod_{m=0}^{n-1}\theta(tq^m; p), \qquad \theta(t; p; q)_0=1,
\label{ePoch}\end{equation}
and use a shorthand notation for the $k$-fold product of elliptic Pochhammer symbols
$$
\theta(t_1, \ldots ,t_k; p; q)_n  :=  \theta(t_1; p; q)_n \cdot  \ldots  \cdot \theta(t_k; p; q)_n.
$$
For $p=0$ we have $\theta(t;0)=1-t$ and the product \eqref{ePoch} reduces to the
standard $q$-Pochhammer symbol
$$
\theta(t; 0; q)_n = (t;q)_n= \prod_{m=0}^{n-1}(1-tq^m).
$$

In this paper we study properties of the elliptic hy\-per\-geo\-met\-ric series. These series  may be considered as a natural generalization of classical objects, such as the generalized hy\-per\-geo\-met\-ric function \cite{aar}
which is defined in terms of the Pochhammer symbol $(a)_n=a(a+1) \ldots (a+n-1)$ and its $k$-fold
product version $(\ell_1, \ldots ,\ell_k)_n=(\ell_1)_n \cdot \ldots \cdot (\ell_k)_n$ as follows
\begin{equation}
_s F_r\left(\begin{matrix} a_0, \ldots , a_{s-1}\\ b_1, \ldots , b_r \end{matrix} \bigg| z\right) = \sum_{n=0}^\infty \frac{(a_0,  \ldots  ,a_{s-1})_n}{(b_1,  \ldots  ,b_r)_n}\frac{z^n}{n!}.
\label{sFr}\end{equation}
This function satisfies the following differential equation of finite order in terms of the operator $\delta=z\frac{d}{dz}$
\begin{equation}
[\delta(\delta+b_1-1) \cdots (\delta+b_r-1)-z(\delta+a_0) \cdots (\delta+a_{s-1})] \cdot _s F_r\left(\begin{matrix} a_0, \ldots , a_{s-1}\\ b_1, \ldots , b_r \end{matrix} \bigg| z\right)=0,
\end{equation}
which is easy to verify in a straightforward way.

A natural generalization of the $_s F_r$-series, which nowadays is also treated as a classical function,
is called the basic, or $q$-hy\-per\-geo\-met\-ric series
\begin{equation}
_s \phi_r\left(\begin{matrix} t_0, \ldots , t_{s-1}\\ w_1, \ldots , w_r \end{matrix}; q \bigg| z\right) = \sum_{n=0}^\infty \frac{(t_0,  \ldots  ,t_{s-1};q)_n}{(w_1,  \ldots  ,w_r;q)_n}
\frac{z^n}{(q;q)_n}.
\label{phi}\end{equation}
This function satisfies the following $q$-difference equation of finite order
\begin{eqnarray}\nonumber &&
[(1-q^\delta)(1-w_1q^{\delta-1})\cdots(1-w_rq^{\delta-1})
\\ && \makebox[2em]{}
-z(1-t_0q^{\delta})\cdots(1-t_{s-1}q^{\delta})]
\, {}_s\phi_r\left(\begin{matrix} t_0, \ldots , t_{s-1}\\ w_1, \ldots , w_r \end{matrix}; q \bigg| z\right)=0,
\label{phi_eq}\end{eqnarray}
where the operator $q^\delta$ acts as a $q$-shift operator on arbitrary functions of $z\in\mathbb{C}$,
$q^\delta f(z):=f(qz)$. The notation $q^{zd/dz}$ can be justified by the action of this
operator on functions admitting the Laurent series expansion,
$$
q^{zd/dz} f(z)=\sum_{k=0}^\infty \frac{(\log q)^k}{k!} \Big(z\frac{d}{dz}\Big)^k \sum_{n\in\mathbb{Z}}c_nz^n
= \sum_{n\in\mathbb{Z}}c_n\sum_{k=0}^\infty \frac{(\log q)^k}{k!} n^k z^n
=f(qz).
$$
For $s=r+1$ our definition of $_{r+1}\phi_r$-series coincides with the one given in \cite{GR},
and for $s\neq r+1$ it differs by the changes $q, t_j, w_k \to q^{-1},  t_j^{-1},  w_k^{-1}$,
and appropriate scaling of the $z$-variable, respectively.

The Pochhammer-Horn type characterization of the described plain and $q$-hy\-per\-geo\-met\-ric
series $\sum_{n\in\ZZ} c_n$ consists in the statement that the ratio of sequential series
coefficients $c_{n+1}/c_n$ is a rational
function of $n$ and $q^n$, respectively. Such series are terminating from below
in a natural way after a special choice of one of the free parameters emerging
from such a definition ($b_0=1$ and $w_0=q$, respectively), which is reflected in the
presence of the $n!$ term in \eqref{sFr} and $(q;q)_n$ term in \eqref{phi}.

The next natural step of generalization of the series of hy\-per\-geo\-met\-ric type, that emerged relatively recently,
is the series characterized by the demand that  $c_{n+1}/c_n$ is either an elliptic function of $n$
with modulus $p=e^{2\pi i\tau}$ or a meromorphic theta function of $n$ with the same modulus \cite{spi:theta1,spi:essays}.
This class of hy\-per\-geo\-met\-ric series degenerates to the basic case for $p=0$ and,
thus, to the plain hy\-per\-geo\-met\-ric series, if  $p=0$, $q \rightarrow 1$.
These infinite series formally satisfy finite difference
equations of infinite order, as described below. Note that a particular terminating elliptic hy\-per\-geo\-met\-ric
series, defining per se an elliptic function of parameters of a specific form, was considered for the first time by
Frenkel and Turaev in \cite{FT}.

\section{Theta hy\-per\-geo\-met\-ric series}

\noindent Consider now the following formal power series \cite{spi:theta1}
\begin{equation}
_s E_r\left({ t_0, \ldots , t_{s-1}\atop w_1, \ldots , w_r}; q, p \bigg| z\right) = \sum_{n=0}^\infty \frac{\theta(t_0,  \ldots  ,t_{s-1}  ; p ; q)_n}{\theta(q, w_1,  \ldots  ,w_r  ; p ; q)_n} z^n,
\label{E}\end{equation}
with the natural restriction on the values of parameters $w_k$ which should not lead to singularites
in the series coefficients, $\theta(w_k;p;q)_n\neq 0$.
For $s=0$ or $r=0$ the parameters $t_j$ or $w_k$ are absent, respectively.\\

\textbf{Lemma 1}. The power series $_s E_r$ satisfies the following formal finite difference equation of
infinite order
\begin{eqnarray} \nonumber &&
(p; p)_\infty ^{s-1-r} \sum_{k_1, \ldots ,k_{r+1} \in \mathbb{Z}}
\prod_{j=1}^r w_j^{k_j} \, (-1)^{k_{r+1}}(-q)^{-\sum_{j=1}^{r} k_j}p^{\sum_{j=1}^{r+1}\binom{k_j}{2}}
f(zq^{\sum_{j=1}^{r+1}k_j})
\\ && \makebox[3em]{}
=z \sum_{\ell_0, \ldots ,\ell_{s-1} \in \mathbb{Z}}
\prod_{j=0}^{s-1}(-t_j)^{\ell_j} \, p^{\sum_{j=0}^{s-1}\binom{\ell_j}{2}}f(zq^{\sum_{j=0}^{s-1}\ell_j}).
\label{eq_ini}\end{eqnarray}
\textit{Proof}. Using properties of the elliptic Pochhammer symbol \eqref{ePoch} one can easily
find the ratio of consecutive coefficients of this series
\begin{align*}
h(n)=\frac{c_{n+1}}{c_n}=
\frac{1}{\theta(q^{n+1};p)}\prod_{m=0}^{s-1} \theta(t_mq^n; p)\prod_{k=1}^r \frac{1}{\theta(w_kq^n;p)}.
\end{align*}
This is a particular meromorphic theta function of the formal variable $n$, since it
satisfies the following quasiperiodicity relations
$$
h\Big(n+\frac{2\pi i}{\log q}\Big)=h(n), \qquad h\Big(n+\frac{\log p}{\log q}\Big)=
(-q^n)^{1+r-s} q\prod_{k=1}^rw_k\prod_{m=0}^{s-1}t_m^{-1}\, h(n),
$$
which is characteristic for theta functions.

Notice now that $\theta(1; p)=0$ and, as follows from \eqref{TPI},
$$
\theta(t q^{\delta}; p) z^n=\frac{1}{(p;p)_\infty}\sum_{k \in \mathbb{Z}} p^{k(k-1)/2}(-t q^{\delta})^k z^n
=\theta(t q^n; p) z^n.
$$
Hence, in terms of theta functions of the $q$-shifting operator $q^\delta$,
we derive the following formal equality
\begin{eqnarray} \nonumber && \makebox[-2em]{}
\left[\theta(q^\delta, w_1 q^{\delta-1},  \ldots  , w_r q^{\delta-1}; p)- z \theta(t_0 q^{\delta},  \ldots  , t_{s-1}q^{\delta}; p)\right] \, {}_s E_r\left(\begin{matrix} t_0, \ldots , t_{s-1}\\ w_1, \ldots , w_r \end{matrix}; q, p \bigg| z\right)
\\  &&
=\sum_{n=1}^\infty \frac{\theta(t_0,  \ldots  ,t_{s-1}  ; p ; q)_n z^n}{\theta(q, w_1,\ldots, w_r;p;q)_{n-1}}
- z\sum_{n=0}^\infty \frac{\theta(t_0,  \ldots  ,t_{s-1}  ; p ; q)_{n+1} z^{n}}
{\theta(q, w_1,\ldots, w_r;p;q)_n}  =0.
\label{elleq}\end{eqnarray}
Expressing theta functions of $q^\delta$ in square brackets as the power series of $q^\delta$
by the use of Jacobi triple product identity \eqref{TPI}, we obtain the desired result.
\hfill$\Box$

\medskip
Equation \eqref{eq_ini} for the simplest series $_0E_0(q, p|z)$ was derived in \cite{spi:bailey}.
For $p=0$  the difference equation of infinite order \eqref{elleq} reduces to the finite order
$q$-difference equation for $q$-hy\-per\-geo\-met\-ric series \eqref{phi_eq}.

We now consider more attentively the coefficients emerging in the infinite bilateral
 sums in \eqref{eq_ini}. For some complex parameters $s_j, \, j=1,\ldots, r,$ one has
 \begin{eqnarray}\nonumber &&
\prod_{j=1}^r \theta(s_j q^{\delta}; p) f(z) = \frac{1}{(p; p)_\infty^r}
\sum_{m_1, \ldots ,m_r \in \mathbb{Z}}\prod_{j=1}^r(-s_j)^{m_j}  p^{\sum_{j=1}^r \binom{m_j}{2}}
f(zq^{\sum_{j=1}^r m_j})
\\ && \makebox[4em]{}
=\frac{1}{(p; p)_\infty^r} \sum_{n\in\ZZ} \Phi_n(s_1, \ldots ,s_r;p) f(zq^n),
\label{coeff_Phi}\end{eqnarray}
where
\begin{equation}
\Phi_n(s_1, \ldots ,s_r;p) :=  \sum_{\substack{m_1, \ldots ,m_r \in \mathbb{Z} \\m_1+ \ldots +m_r=n}} (-s_1)^{m_1} \ldots (-s_r)^{m_r} p^{\binom{m_1}{2}+ \ldots +\binom{m_r}{2}}.
\label{Phi_n}\end{equation}
As one may notice, $\Phi_n(s_1, \ldots ,s_r;p)$ is a certain theta function, which depends on $n$ somehow,
and $\Phi_0(s_1, \ldots ,s_r;p)$ is a theta function on the root system $A_{r-1}$.\\

\textbf{Lemma 2}. The following three identities hold true
\begin{align}
&\Phi_n(s_1, \ldots ,s_r;p) = p^{\binom{n}{2}}(-s_j)^n \Phi_0(s_1, \ldots ,s_j p^n, \ldots ,s_r;p),
\quad j=1,\ldots, r,
\\
&\Phi_{n+r}(s_1, \ldots ,s_r;p) =(-1)^rs_1 \ldots s_r p^n \Phi_{n}(s_1, \ldots ,s_r;p),
\label{Phinr} \\
&\Phi_0(s_1, \ldots ,s_r; p) = \frac{1}{r}\frac{(p;p)_\infty^r}{(p^r; p^r)_\infty}
\frac{\sum_{m=0}^{r-1} \theta(s_1 z \zeta^m; p)\cdot  \ldots  \cdot\theta(s_r z \zeta^m; p)}
{\theta(-(-z)^rs_1\cdot  \ldots  \cdot s_r; p^r)},
\label{Phi_sum}
\end{align}
where in the third equality $z$ is an arbitrary non-zero variable and $\zeta$ is any primitive root of
unity, $\zeta^r=1$ (e.g., $\zeta=e^{2\pi i/r}$).\\

\noindent \textit{Proof}. For the first relation, notice the following transformation
\begin{align*}
\Phi_n(s_1, \ldots ,s_r;p)&=\sum_{\substack{m_1, \ldots ,m_r \in \mathbb{Z} \\m_1+ \ldots +m_r=n}} (-s_1)^{m_1} \ldots (-s_r)^{m_r} p^{\binom{m_1}{2}+ \ldots +\binom{m_r}{2}}\\
& \makebox[-5em]{}
=\sum_{\substack{m_1, \ldots ,m_r \in \mathbb{Z} \\m_1+ \ldots +m_r=0}} (-s_1)^{m_1} \ldots (-s_j)^{m_j+n} \ldots (-s_r)^{m_r} p^{\binom{m_1}{2}+ \ldots +\binom{m_j+n}{2}+ \ldots +\binom{m_r}{2}}
\\
&= p^{\binom{n}{2}}(-s_j)^n \Phi_0(s_1, \ldots ,s_j p^n, \ldots ,s_r;p).
\end{align*}
Similar computation then yields equality \eqref{Phinr}
\begin{align*}
\Phi_{n+r}(s_1, \ldots ,s_r;p)&=\sum_{\substack{m_1, \ldots ,m_r \in \mathbb{Z} \\m_1+ \ldots +m_r=n+r}} (-s_1)^{m_1} \ldots (-s_r)^{m_r} p^{\binom{m_1}{2}+ \ldots +\binom{m_r}{2}}\\
&=\sum_{\substack{m_1, \ldots ,m_r \in \mathbb{Z} \\m_1+ \ldots +m_r=n}} (-s_1)^{m_1+1} \ldots (-s_r)^{m_r+1} p^{\binom{m_1+1}{2}+ \ldots +\binom{m_r+1}{2}}\\
&= (-1)^rs_1 \ldots s_r p^n \Phi_n(s_1, \ldots ,s_r;p).
\end{align*}
Notice that the latter identity implies the relation
\begin{equation}
\Phi_{nr}(s_1, \ldots ,s_r; p)=p^{r\binom{n}{2}}(-1)^{nr}(s_1 \ldots s_r)^n\Phi_0(s_1, \ldots ,s_r;p).
\label{Phi_nr}\end{equation}

For the last part \eqref{Phi_sum}, consider now the  expansion similar to \eqref{coeff_Phi}
\begin{align*}
\theta(s_1 z; p)\cdot  \ldots  \cdot\theta(s_r z; p) &
=\frac{1}{(p;p)_\infty^r}\sum_{n=-\infty}^{\infty} z^n \Phi_n(s_1, \ldots ,s_r;p).
\end{align*}
After using the property of primitive $r$-th roots of unity $\sum_{k=0}^{r-1}\zeta^{kn}=0$ for $n\neq \ell r$,
a Fourier-type transformation then gives
\begin{align*}
& \frac{1}{r}\sum_{m=0}^{r-1} \theta(s_1 z \zeta^m; p)\cdot  \ldots  \cdot\theta(s_r z \zeta^m; p) = \frac{1}{(p;p)_\infty^r}\sum_{n=-\infty}^{\infty} z^{nr}\Phi_{nr}(s_1, \ldots ,s_r; p)
\\ & \makebox[4em]{}
=\frac{(p^r;p^r)_\infty}{(p;p)_\infty^r} \theta(-(-z)^r s_1 \ldots s_r; p^r)\Phi_{0}(s_1, \ldots ,s_r; p),
\end{align*}
where we used relation \eqref{Phi_nr}.
The explicit expression for function $\Phi_0(s_1, \ldots ,s_r; p)$ now obviously follows. \hfill$\Box$

Finally, we come to the following statement.\\

\textbf{Theorem 3}. The theta hy\-per\-geo\-met\-ric series
\eqref{E} is a formal power series solution of the infinite order $q$-difference equation
\begin{equation}
\sum_{n\in\ZZ}\Big( \Phi_n(1,q^{-1}w_1, \ldots ,q^{-1}w_r;p)- z\Phi_n(t_0, \ldots ,t_{s-1};p)\Big) f(zq^n)=0,
\label{final_eq}\end{equation}
with the coefficients $\Phi_n$ defined in \eqref{Phi_n}.\\

Equation \eqref{final_eq} is obtained simply after applying the operator relation \eqref{coeff_Phi}
to both operator theta function products standing in the square brackets in the first line
of equation \eqref{elleq}.

For example, for $r=2$ we can find an explicit expression for the function $\Phi_0(s_1, \ldots ,s_r;p)$ that may be obtained from the representation \eqref{Phi_sum} after substituting $z=s_2^{-1}$,
$$
\Phi_0(s_1,s_2;p)=\sum_{m\in\ZZ} \Big(\frac{s_1}{s_2}\Big)^m p^{m^2}
= (p^2;p^2)_\infty \theta\left(-p\frac{s_1}{s_2};p^2\right).
$$
This means that the formal theta hy\-per\-geo\-met\-ric series ${}_{2}E_1({a,b\atop c}\big|z)$
satisfies the following explicit compact equation
\begin{equation}
\sum_{n\in\ZZ} (-1)^n p^{\binom{n}{2}} \left( \theta\left(-\frac{q}{c}p^{n+1};p^2\right)
- za^n \theta\left(-\frac{a}{b}p^{n+1};p^2\right)  \right)  {}_2 E_1\left({ a,b\atop c}; q, p \bigg| q^n z\right)=0,
\label{2E1_eq}\end{equation}
which can be considered as a direct theta functional analogue of the
hy\-per\-geo\-met\-ric equation for the Euler-Gauss $_2F_1$-function. It can take different
forms obtained after choosing different parametrization of the $\Phi_n$-coefficients.
In particular, one can permute the parameters $a$ and $b$.
Note that in all cases the $q$-shifted functions $f(q^n z)$ enter the equation
only with a coefficient linear in the independent variable $z$. Evidently, any solution
of \eqref{2E1_eq} can be multiplied by an arbitrary elliptic function
$\varphi(z)=\varphi(qz)$ and it remains a solution of the same equation.

Splitting the summation index $n$ in \eqref{2E1_eq} to even and odd values
and clearing $p^n$-factors in the theta function arguments, we come to the equation
\begin{eqnarray} \nonumber &&  \makebox[-1em]{}
\sum_{k\in\ZZ}\bigg\{ p^{k(k-1)}\bigg[ \theta\left(-\frac{pq}{c};p^2\right)\left(\frac{c}{q}\right)^k
-z\theta\left(-\frac{pa}{b};p^2\right)(ab)^k  \bigg] f(q^{2k}z)
\\ &&
-p^{k^2}\bigg[ \theta\left(-\frac{q}{c};p^2\right)\left(\frac{c}{q}\right)^{k+1}
-zb\theta\left(-\frac{a}{b};p^2\right)(ab)^k  \bigg] f(q^{2k+1}z)\bigg\}=0.
\label{2E1_eqSPLIT}\end{eqnarray}
For fixed parameters $a,b,c,d,$ and $q$, only the terms with the summation
index $k=0,1$ give nontrivial contributions in the $p\to 0$ limit,
provided the function $f(z)$ depends on $p$ in a non-singular way. This yields
the $q$-hy\-per\-geo\-met\-ric equation for the function $_2\varphi_1$ \eqref{phi_eq}.

Supposedly, a finite-difference equation of infinite order should have rich space of solutions
with many constants of integration. At the moment we cannot characterize the whole
space of solutions of equation \eqref{2E1_eqSPLIT}. However, it is possible to find solutions with
the isolated singularity at the point $z=0$. Indeed, any such function can be expanded around
$z=0$ into the Laurent series $f(z)=\sum_{n\in\ZZ}c_n z^n$. Substituting this
expansion into  \eqref{2E1_eqSPLIT} one finds a recurrence relation of the first order
for the coefficients $c_n$ which can be solved explicitly in terms of the elliptic Pochhammer
symbols. For generic values of parameters the coefficient $d_0$ in the relation
$c_{-1}=d_0 c_0$ appears to be equal to zero, $d_0=0$, which means that all $c_n=0$ for $n<0$.
Therefore the whole series appears to be analytical at the point $z=0$ (provided it
converges) and it coincides with the derived above solution
$f(z)={}_2 E_1\big({ a,b\atop c}; q, p \big| z\big)$. Multiplying
this solution by a nontrivial elliptic function $\varphi(qz)=\varphi(z)$, we obtain a single-valued
solution of \eqref{2E1_eq} which is not analytical at $z=0$. Other solutions can contain
the parts which are not single valued functions of $z$.

Some particular solutions for toy equations of infinite order are described below.
Consider eigenfunctions of the simplest difference operator of infinite order $\theta(aq^\delta; p)$.
It should give a feeling on the possible structure of solutions of more complicated $q$-difference
equations of infinite order. Using the infinite product representation instead of the
infinite theta series form, we have
$$
\theta(aq^\delta; p) f(z)=\prod_{n=0}^\infty (1-aq^\delta p^n) (1-pa^{-1}q^{-\delta}p^n)f(z)= \lambda f(z).
$$
A single-valued function $f(z)=z^N\theta(bz;q)/\theta(cz;q), \, N\in\ZZ,$ satisfies the equation $f(qz)=q^Ncb^{-1}f(z)$.
Therefore we have an eigenfunction
\begin{equation}
f(z)=z^N\prod_{j=1}^m\frac{\theta(b_jz;q)}{\theta(c_jz;q)}, \qquad
\lambda=\theta\big(aq^N\prod_{j=1}^m \frac{c_j}{b_j};p\big),
\label{eigenf}\end{equation}
for arbitrary sets of parameters $b_j, c_j\in \CC^\times$.
Whenever the condition $\prod_{j=1}^m(b_j/c_j)=aq^Np^n,\, n\in\ZZ,$ is satisfied,
one gets a kernel function. As a result, any linear combination of such functions with
different parameters $N, b_j, c_j$ satisfying the same constraint also lies in the kernel of the taken operator.

There are also solutions which are not single-valued  in $z$ similar to the general solution
of the ordinary hy\-per\-geo\-met\-ric equation for the $_2F_1$
Euler-Gauss hy\-per\-geo\-met\-ric function. Indeed, consider the function of the form
$z^d$, which is not single valued for $d\neq \ZZ$. One has the formal eigenvalue
problem $\theta(aq^\delta; p) z^d= \theta(aq^d; p) z^b$, so that whenever $aq^d=p^n$
(or $d=-\log ap^n/\log q$), $n\in\ZZ$, one obtains a kernel function. Their linear combination with different
values of $n$ weighted by $q$-elliptic coefficients yields the following infinite Dirichlet-type  series
\begin{equation}
f(z)=z^{-\frac{\log a}{\log q}}\sum_{k\in \mathbb{Z}} h_k(z) z^{-k \frac{\log p}{\log q}},
\quad h_k(qz)=h_k(z),
\label{kernel}\end{equation}
which will be a kernel function of the operator $\theta(aq^\delta; p)$, provided the series converges.
One can replace also in the single-valued function \eqref{eigenf} the factor $z^N$ by $z^d$,
impose the constraint $\prod_{j=1}^m(b_j/c_j)=aq^dp^n,\, n\in\ZZ,$ and take linear
combinations of such functions for different $b_j, c_j, d$ and $n$ to obtain even more complicated kernel functions.
 Suppose that the operator $\theta(aq^\delta; p)$
acts in some Hilbert space of functions and one is interested in its inversion. If at least one of the
described kernel functions  belongs to the Hilbert space of interest, then the inversion is not possible.

Consider now a more complicated operator. Take $w, q, p\in \CC^\times, \; K\in\ZZ$,
set $t:=pq^Kw^{-K^2}$ and impose the constraints $|w|, |p|, |t|<1$. Then define the operator
$$
\hat{\mathcal{O}}=\sum_{n\in \mathbb{Z}} w^{\binom{n}{2}}\theta(w^{nK}z;p) q^{n\delta}.
$$
Expanding $\theta(w^{nK}z;p)$ into the Laurent series in $z$ and changing the order of
summations in the emerging double sum, we obtain
$$
\hat{\mathcal{O}}=\frac{(w;w)_\infty (t;t)_\infty}{(p;p)_\infty}
\theta(w^{-\binom{K}{2}}zq^{-K\delta};t)\theta(-q^\delta;w),
$$
where we used the relation $z^mq^{-Km\delta}=q^{K\binom{m}{2}}(zq^{-K\delta})^m$.
As we see, the operator of the form $\theta(az q^{K\delta}; t)\theta(bq^\delta; w)$
determines the kernel of the $\hat{\mathcal{O}}$-operator and the functions described above
provide explicit examples of that.

\section{Convergence of the elliptic hy\-per\-geo\-met\-ric series}

We now return to examination of the $_s E_r$-series, in particular, we raise the question, when this formal power series converges. One may readily notice that substitution $t_i = q^{-N}$ for some fixed $i$ and
$N \in \mathbb{Z}_{\geq 0}$ leads to the termination of the series
$$
_s E_r\left(\begin{matrix} q^{-N}, t_1, \ldots , t_{s-1}\\ w_1, \ldots , w_r \end{matrix}; q, p \bigg| z\right) =
\sum_{n=0}^{N} \frac{\theta(q^{-N}, t_1, \ldots ,t_{s-1}  ; p ; q)_n}{\theta(q, w_1,  \ldots  ,w_r  ; p ; q)_n} z^n,
$$
and the resulting sum converges for a trivial reason.

The most interesting case corresponds to the elliptic hy\-per\-geo\-met\-ric series.
Set $s=r+1$ and impose special balancing constraint
\begin{equation}
\prod_{k=0}^{r} t_k = q\prod_{k=1}^{r} w_k.
\label{balance}\end{equation}
Then we see that the ratio of sequential series terms
$$
h(n)=\frac{c_{n+1}}{c_n}
=z\frac{\theta(q^n t_0, q^n t_1,  \ldots  , q^n t_{r}  ; p)}{\theta(q^{n+1}, q^n w_1,  \ldots  , q^n w_r  ; p)}
$$
becomes an elliptic function of the formal variable $n$. Namely, considering $n$ as a complex variable,
$h(n)$ becomes a meromorphic function of $n$ with double periodicity
\begin{equation}
h\Big(n+\frac{2\pi i}{\log q}\Big)=h\Big(n+\frac{\log p}{\log q}\Big)=h(n),
\label{ellhn}\end{equation}
i.e. it is an elliptic function. The latter fact gives the name to the series. Equalities \eqref{ellhn}
follow from the fact that $h(n)$ depends only on the function $q^n$, which does not change under
the first shift of the $n$-variable and it becomes $pq^n$ after the second period shift.
In the latter case the quasiperiodicity of theta functions \eqref{trafo} brings in constant multipliers
which cancel due to the balancing condition \eqref{balance}.

In order to remove singularity of $h_n$ at $n=-1$ we take now special value $t_0=q$ and consider the series
$$
_{r+1} E_r\left(\begin{matrix} q, t_1, \ldots , t_{r}\\ w_1, \ldots , w_r \end{matrix}; q, p \bigg| z\right) = \sum_{n=0}^\infty \frac{\theta(t_1,  \ldots  ,t_{r}  ; p ; q)_n}{\theta(w_1,  \ldots  ,w_r  ; p ; q)_n} z^n,
\quad \prod_{k=1}^{r} t_k = \prod_{k=1}^{r} w_k.
$$
The balancing condition implies that the function
\begin{equation}
H(u):=\frac{\theta(u t_1,  \ldots  , u t_{r}  ; p)}{\theta(u w_1,  \ldots  , u w_r  ; p)}
\label{Hu}\end{equation}
is a meromorphic function of $u\in\mathbb{C}^\times$ invariant under the transformation $u\to pu$
(or $p$-elliptic function),
$H(pu)=H(u)$. Zeros of this function are located at the points $u=t_k^{-1} p^\mathbb{Z}$  and its poles sit at $u=w_k^{-1}p^{\mathbb{Z}}$ for $k=1, \ldots ,r$.

In the exponential form, we denote
$$
p=e^{2\pi i\tau}, \quad \mathrm{Im}\,\tau>0,
$$
and then the function $H(e^{2\pi i x })$ is an elliptic function of $x\in \mathbb{C}$
with the fundamental periods $\{1, \tau\}$.
We now carefully choose $q$. Consider any pair of integers $N, M$, such that the line from $x=0$ to $x=N+M\tau$
does not contain any poles of the function $H(e^{2\pi i x })$.
For such a pair of integers and any real number $\chi \in \mathbb{R}$ we set
\begin{equation}
q=e^{2\pi i \chi(N+M\tau)}, \qquad (N,M)\neq (0,0).
\label{q_choice}\end{equation}
This assumes that for fixed $\chi  \in \mathbb{R}$ and integers
$N$ and $M$ the parameters satisfy the constraints
\begin{equation}
t_j,\, w_j\neq p^k  e^{2\pi i\chi  (N+M\tau)l},
\quad k\in{\mathbb{Z}},\; l\in \mathbb{Z}_{\leqslant 0}, \quad j=1,\ldots, r.
\label{tw}\end{equation}
Then, the series coefficients $c_n/z^n$ neither vanish (i.e. the series does not terminate), nor
blow up. However, if some of the parameters $w_j$ lie on the line from $x=0$ to $x=N+M\tau$ at the
transcendental points, then in the limit $n\to\infty$ some subsequences of the points $q^np^{\ZZ}$ will approach
these parameters arbitrarily close, there emerges a singularity and the estimate of the values of $|H(q^n)|$
becomes delicate.
Therefore, in the simplest case one can assume that the constraints on parameters $w_j$
\eqref{tw} are true for all real values of $\chi $, not just for a fixed number. Then the absolute
values $|H(q^n)|$ are bounded from above by some positive constant $S$. As a result the series
converges for $|z| < S^{-1}$.

Suppose that $\chi $ is a rational number, $\chi =a/b,\, (a,b)=1$, i.e. $q^b=p^{Ma}$. Then one has the
periodicity $H(q^{n+b})=H(q^n)$ and the series of interest can be summed to a closed form
\begin{eqnarray}\nonumber &&
\sum_{n=0}^\infty \frac{\theta(t_1,  \ldots  ,t_{r}  ; p ; q)_n}{\theta(w_1,  \ldots  ,w_r  ; p ; q)_n} z^n
=\sum_{j=0}^\infty\sum_{l=0}^{b-1}\prod_{k=0}^{jb+l-1}H(q^k)z^{jb+l}
\\  && \makebox[-1em]{}
=\sum_{j=0}^\infty R^j z^{jb}\sum_{l=0}^{b-1}\prod_{k=jb}^{jb+l-1}H(q^k)z^{l}
=\frac{1+\sum_{l=1}^{b-1}\prod_{k=0}^{l-1} H(q^{k}) z^l}{1-Rz^b},
\label{roots}\end{eqnarray}
where we have denoted $R:=\prod_{k=0}^{b-1} H(q^k)$. Evidently, the radius of convergence
in this case is $r_{c}=|R|^{-1/b}$. Since for nontrivial $R$-function $|R|$ cannot be equal to 1 for
all values of parameters, whenever $|R|>1$ we can flip values of $t_k$ and $w_k$ and obtain $r_c>1$,
i.e. there always exist series with such a radius of convergence. Note that formula \eqref{roots}
defines a meromorphic function of parameters $t_k, w_k, q$ as a rational combination of theta
functions, and in general it is neither elliptic, nor even quasiperiodic in these variables.

Consider now the case when $\chi $ is a fixed irrational number, with additional constraints
on the values of $t_j$, similar to those of $w_j$, so that not only the poles of the function $H(u)$ do not lie on the considered interval, but the zeros $t_j^{-1}p^{\mathbb{Z}}$ are also assumed not to lie at transcendental points.
In this case the following explicit bound on the radius of convergence may be obtained.
By the Weyl equidistribution theorem for any irrational $\chi $ the logarithms of the numbers $q^k = e^{2\pi i k\chi (N+M\tau)}$ are uniformly distributed on the line from $0$ to $N+M\tau$. Therefore, the following limit exists
\begin{equation}
\lim_{n \to \infty} \frac{1}{n} \sum_{k=0}^{n-1} \log |H(q^k)| = \int_{0}^{1} \log |H(e^{2\pi i x (N+M\tau)})|d x.
\label{sum_int}\end{equation}
As a result, the leading asymptotics of our series coefficients has the plain power form
$$
|c_n| =\prod_{k=0}^{n-1}|H(q^k)|\, |z|^n
=\big|z\exp \big(\frac{1}{n} \sum_{k=0}^{n-1} \log |H(q^k)|\big)\big|^n
\stackreb{=}{n\to\infty}  \big|\frac{z}{r_c}\big|^n,
$$
where
$$
r_{c}:
= \exp\left(-\int_{0}^{1} \log |H(e^{2\pi i x (N+M\tau)})|d x \right).
$$
By the ratio test, we see that our series converges for $|z|<r_c$, i.e. $r_c$ is the radius of
convergence.

The convergence problem for a particular
$q$-hy\-per\-geo\-met\-ric series when $q = e^{2\pi i \chi }$ for irrational $\chi $ was considered
for the first time by Hardy and Littlewood \cite{HL}, where a simpler integral emerged and it
was explicitly evaluated giving exact value for the radius of convergence.
However, in \cite{HL} an analogue of $H(u)$ had a
singularity and uniform distribution worked not for all irrational $\chi $.
For more results of this type for $q$-hy\-per\-geo\-met\-ric objects, see \cite{Lub,Petr}
and the next section. In particular,  Petruska has shown in \cite{Petr} that it is essential
for $t_i, w_j$ not to lie in the interval of integration, which we assumed.

One can pose a natural question --- if there is any way to compute the arising integral?
Consider the following analogue of the above integral, where we assume $|p|<1$, or $\mathrm{Im}\,\,\tau>0$,
\begin{equation}
F_{N,M}(t):=\int_{0}^{1} \log |\theta(te^{2\pi i x  (N+M\tau)}; p)|d x.
\label{FNMdef}\end{equation}
Singularities of the integrand function may only be of logarithmic type,
hence the integral converges. We first examine dependence of this integral
on the parameters $N$ and $M$. For that purpose consider their greatest common
divisor $D=\gcd(N,M)$, $D>0$. Then
\begin{equation}
F_{N,M}(t)=F_{\frac{N}{D}, \frac{M}{D}}(t)+\pi \mathrm{Im}\,\tau \frac{M(D-1)(2M(D+1)-3D)}{6D^2}-\log|t|\frac{M(D-1)}{2D}.
\label{FNM/D}\end{equation}
This equality is obtained in a consequence of several elementary transformations. Starting from
\begin{align*}
&D\int_{0}^{1} \log |\theta (te^{2\pi i x  (N+M\tau)}; p)|d x =\int_{0}^{D} \log|\theta (te^{2\pi i  x (\frac{N}{D}+\frac{M}{D}\tau)};p)|d x
\\ &   =\sum_{k=0}^{D-1} \int_{0}^{1}\log|\theta (te^{2\pi i ( x +k)(\frac{N}{D}+\frac{M}{D}\tau)};p)|d x
=\sum_{k=0}^{D-1} \int_{0}^{1}\log|\theta (p^{Mk/D}te^{2\pi i x (\frac{N}{D}+\frac{M}{D}\tau)};p)|d x,
\end{align*}
applying the quasiperiodicity relation \eqref{TPI} and taking the logarithm of the module
of the quasiperiodicity multiplier, we obtain
\begin{align*}
F_{N,M}(t)=F_{\frac{N}{D}, \frac{M}{D}}(t)+\frac{1}{D}\sum_{k=0}^{D-1} \left(\frac{M^2k(k+1)}{D^2}\pi\mathrm{Im}\,\tau-\frac{Mk}{D}\pi \mathrm{Im}\,\tau-\frac{Mk}{D}\log|t|\right).
\end{align*}
Computing the sums over $k$, we come to the claimed result.

Next, we describe a different representation of $F_{N,M}(t)$ based on the infinite product form
of the theta-function. We restrict our attention to the case $M>0$ and write
\begin{align*} & \makebox[-1em]{}
\int_{0}^{1} \log |\theta (te^{2\pi i x  (N+M\tau)}; p)|d x
=\int_{0}^{1} \sum_{k=0}^{M-1}  \log |\theta (tp^k e^{2\pi i x  (N+M\tau)}; p^M)|d x
\\ &
=\sum_{k=0}^{M-1}  \int_{0}^{1} \sum_{n=0}^\infty\log |(1-tp^k e^{2\pi i( x +n)(N+M\tau)})(1-t^{-1}p^{-k}
e^{2\pi i(n+1- x )(N+M\tau)})|d x  \\
&   =\sum_{k=0}^{M-1}  \sum_{n=0}^\infty \int_{n}^{n+1}\log |(1-tp^k e^{2\pi i x (N+M\tau)})(1-t^{-1}p^{-k} e^{2\pi i x (N+M\tau)})|d x \\
&   =\sum_{k=0}^{M-1}  \int_{0}^\infty \log |(1-tp^k e^{2\pi i  x (N+M\tau)})
(1-t^{-1}p^{-k} e^{2\pi i x (N+M\tau)})|d x.
\end{align*}
Hence our problem has reduced to an explicit computation of the following integral for arbitrary $t\in \mathbb{C}^\times$:
$$
I_{N,M}:=\int_{0}^\infty \log|1-te^{2\pi i x(N+M\tau)}|dx.
$$
We notice that the case $M<0$ is equivalent to this one because of the transformation $ x  \leftrightarrow 1- x $ in the initial integral. We return to the case $M=0$ later. For now consider the following lemma.

\textbf{Lemma 4}. Suppose $M >0$ and define
\begin{equation}
\mu := \frac{t}{|t|}e^{i (N+M\mathrm{Re}\, \tau)\frac{\log|t|}{M\mathrm{Im}\,\tau}},
\quad |\mu|=1.
\label{mu}\end{equation}
Then the following holds:
\begin{align*}
&\text{if} \quad |t|\geqslant 1, \quad\text{then} \quad
I_{N,M}= \frac{\log^2|t|}{4\pi M\mathrm{Im}\,\tau} - \mathrm{Re} \frac{\mathrm{Li}_2(t^{-1})}{2\pi i (N+M\tau)}-\frac{M\mathrm{Im}\,\tau \mathrm{Re}(\mathrm{Li}_2(\mu))}{\pi|N+M\tau|^2},
\\
&\text{if}\quad|t|\leqslant 1, \quad\text{then} \quad
I_{N,M}= \mathrm{Re} \frac{\mathrm{Li}_2(t)}{2\pi i (N+M\tau)},
\end{align*}
where $\mathrm{Li}_2(x)$ is the classical dilogarithm function
$$
\mathrm{Li}_2(x)=\sum_{n=1}^\infty \frac{x^n}{n^2}, \quad x\in\CC, \; |x|\leqslant 1.
$$
\textit{Proof}. Consider the first case, i.e. $|t|\geqslant 1$. Then the integral may be divided into two pieces,
$$
I_{N,M} = \int_{0}^{\frac{\log|t|}{2\pi M\mathrm{Im}\,\tau}} \log|1-te^{2\pi i x(N+M\tau)}|dx + \int_{\frac{\log|t|}{2\pi M\mathrm{Im}\,\tau}}^\infty \log|1-te^{2\pi i x(N+M\tau)}|dx.
$$
For the first term we notice that $0\leqslant x \leqslant \frac{\log|t|}{2\pi M\mathrm{Im}\,\tau}$. Therefore $\log|t|-2\pi Mx\mathrm{Im}\,\tau \geqslant 0$ or, in other words, $|te^{2\pi i x(N+M\tau)}|=|t|e^{-2\pi Mx\mathrm{Im}\,\tau}\geqslant 1$. Hence we can write
\begin{align*}
\log|1-te^{2\pi i x(N+M\tau)}|=\log|t|-2\pi Mx\mathrm{Im}\,\tau+ \log|1-t^{-1}e^{-2\pi i x(N+M\tau)}|
\end{align*}
and expand the logarithm on the right-hand side into the convergent Taylor series. After taking the termwise
integrals over $x$, this yields
\begin{eqnarray} \nonumber
&\int_{0}^{\frac{\log|t|}{2\pi M\mathrm{Im}\,\tau}} \log|1-te^{2\pi i x(N+M\tau)}|dx
= \frac{\log^2|t|}{4\pi M\mathrm{Im}\,\tau}-\mathrm{Re}\, \sum_{n=1}^\infty \frac{t^{-n}}{n}
 \frac{1-e^{-in(N+M\tau)\frac{\log|t|}{M\mathrm{Im}\,\tau}}}{2\pi in(N+M\tau)}
\\ &
\qquad \quad= \frac{\log^2|t|}{4\pi M\mathrm{Im}\,\tau}-\mathrm{Re}  \frac{\mathrm{Li}_2(t^{-1})}{2\pi i(N+M\tau)}+\mathrm{Re} \frac{\mathrm{Li}_2(\bar \mu)
}{2\pi i(N+M\tau)},
\label{firstint}\end{eqnarray}
where $\bar\mu$ is the complex conjugate of the parameter defined in the Lemma statement.
The last summand is defined correctly, since the argument lies on the unit circle.

Similarly, for the second summand $|te^{2\pi i x(N+M\tau)}|=|t|e^{-2\pi Mx\mathrm{Im}\,\tau}\leqslant 1$, so
that we can expand the logarithm into Taylor series and integrate termwise, which yields
\begin{align*}
\int_{\frac{\log|t|}{2\pi M\mathrm{Im}\,\tau}}^\infty \log|1-te^{2\pi i x(N+M\tau)}|dx
= \mathrm{Re}\, \sum_{n=1}^\infty \frac{t^n}{n}\frac{e^{i n(N+M\tau)\frac{\log|t|}{M\mathrm{Im}\,\tau}}}{2\pi i n(N+M\tau)}
= \mathrm{Re} \frac{\mathrm{Li}_2(\mu)
}{2\pi i (N+M\tau)}.
\end{align*}
Thus, the full integral is indeed equal to
$$
\int_{0}^\infty \log|1-te^{2\pi i x(N+M\tau)}|dx = \frac{\log^2|t|}{4\pi M\mathrm{Im}\,\tau}-\mathrm{Re}  \frac{\mathrm{Li}_2(t^{-1})}{2\pi i(N+M\tau)}+\mathrm{Re} \frac{\mathrm{Li}_2(\bar{\mu})
+\mathrm{Li}_2(\mu)}{2\pi i(N+M\tau)}.
$$

Consider the second case, i.e. $|t|\leqslant 1$. Then $\frac{\log|t|}{2\pi M\mathrm{Im}\,\tau} \leqslant x$ for any
$0\leqslant x$ by default. In other words $|te^{2\pi i x(N+M\tau)}|=|t|e^{-2\pi Mx\mathrm{Im}\,\tau}\leqslant 1$.
Hence the following holds
\begin{equation}\nonumber
\int_0^\infty \log|1-te^{2\pi i x(N+M\tau)}|dx = -\mathrm{Re}\, \sum_{n=1}^\infty \frac{t^n}{n}
\frac{e^{2\pi i nx(N+M\tau)}}{2\pi i n(N+M\tau)}\bigg|_0^\infty
= \mathrm{Re} \frac{\mathrm{Li}_2(t)}{2\pi i(N+M\tau)}.
\end{equation}
\hfill$\Box$

\textit{Remark}.
By definition we have $|\mu|=1$ and $\mathrm{arg} \, \mu \in [0; 2\pi]$. Therefore, the following holds true
\begin{equation}
\mathrm{Re}(\mathrm{Li}_2(\mu))=\pi^2B_2\left(\frac{\arg \mu}{2\pi}\right )= \frac{(\mathrm{arg}\, \mu)^2}{4}-\frac{\pi \mathrm{arg} \, \mu}{2}+\frac{\pi^2}{6}
\label{bern}\end{equation}
for
$$
 \mathrm{arg} \, \mu= \frac{\mathrm{Re} ((N+M\tau)\log \bar{t})}{M\mathrm{Im}\,\tau} \mod 2\pi,
$$
where we use the standard Fourier expansion of the second Bernoulli polynomial:
\begin{eqnarray*} &&
a_n:=\int_0^1B_2(x)e^{-2\pi i nx}dx=
\left\{
\begin{array}{cl}
0, &  n=0  \\ [0.5em]   \displaystyle
\frac{1}{2\pi^2 n^2}, & n\in\ZZ /\{0\}
\end{array}
\right. ,
\\ &&  \makebox[3em]{}
B_2(\{x\})=\sum_{n\in\ZZ}a_ne^{2\pi i nx}=\{x\}^2-\{x\}+\frac{1}{6}.
\end{eqnarray*}

The latter lemma implies that for $M\neq 0$ and any integer $N$ the integral $F_{N,M}(t)$ may always be computed. We have
\begin{eqnarray} \nonumber &&
\int_{0}^\infty \log|(1-te^{2\pi i x(N+M\tau)})(1-t^{-1}e^{2\pi i x(N+M\tau)})|dx
\\ &&  \makebox[2em]{}
=\frac{\log^2|t|}{4\pi M\mathrm{Im}\,\tau}-\frac{M\mathrm{Im}\,\tau \,\mathrm{Re}(\mathrm{Li}_2(\mu))}{\pi|N+M\tau|^2},
\quad |t|\geqslant 1.
\label{integral}\end{eqnarray}
Since $\mu(t)$, as a function of $t$, satisfies relation $\mu(t^{-1})=\overline{\mu(t)}$,
for $|t| \leqslant 1$ we have exactly the same answer.
As a result, from the relation $\mu(tp^k)=\mu(t)e^{-2\pi i k\frac{N}{M}}$, for any
$t \in \mathbb{C}^\times$  we have
$$
\int_{0}^{1} \log |\theta(te^{2\pi i x  (N+M\tau)}; p)|d x  =\sum_{k=0}^{M-1} \left(\frac{\log^2|tp^k|}{4\pi M\mathrm{Im}\,\tau}-\frac{M\mathrm{Im}\,\tau \mathrm{Re}(\mathrm{Li}_2(\mu(tp^k)))}{\pi|N+M\tau|^2}\right).
$$

Consider again the greatest common divisor of $N$ and $M$, $D=\gcd(N,M)$, where we assume as usual that $\gcd(N,0)=N$ and $\gcd(0,M)=M$. Then for the second summand we have the following computation
\begin{align*}
\sum_{k=0}^{M-1} & \mathrm{Li}_2(\mu(tp^k)) =  \sum_{k=0}^{M-1} \sum_{n=1}^\infty \frac{e^{in\,\mathrm{arg} \, \mu(tp^k)}}{n^2}=
\sum_{k=0}^{M-1} \sum_{n=1}^\infty \frac{e^{in\,\mathrm{arg} \, \mu-2\pi i kn\frac{N}{M}}}{n^2}
\\
&= \sum_{k=0}^{M-1} \sum_{\substack{n=1 \\ \frac{nN}{M} \in \mathbb{Z}}}^\infty
\frac{e^{in\,\mathrm{arg} \, \mu-2\pi i kn\frac{N}{M}}}{n^2}=  \sum_{k=0}^{M-1} \sum_{r=1}^\infty \frac{e^{ir\frac{M}{D}\,\mathrm{arg} \, \mu-2\pi i kr\frac{M}{D}\frac{N}{M}}}{(r\frac{M}{D})^2}
=\frac{D^2}{M} \mathrm{Li}_2(\mu^{\frac{M}{D}}).
\end{align*}
And the first summand is clearly equal to the following expression
$$
\frac{\log^2|t|}{4\pi \mathrm{Im}\,\tau}+\frac{\log|t|\log|p|}{4\pi \mathrm{Im}\,\tau}(M-1)
+\frac{\log^2|p|}{4\pi \mathrm{Im}\,\tau}\frac{(M-1)(2M-1)}{6}.
$$
Combining them together we obtain the following result for function $F_{N,M}(t)$ \eqref{FNMdef} with $M>0$:
\begin{eqnarray}\nonumber &&
F_{N,M}(t) = \frac{\log^2|t|}{4\pi\mathrm{Im}\,\tau}-\frac{(M-1)\log|t|}{2}
+\frac{(M-1)(2M-1)\pi\mathrm{Im}\,\tau}{6}
\\ && \makebox[4em]{}
-\frac{D^2\mathrm{Im}\,\tau}{\pi|N+M\tau|^2}
\mathrm{Re}(\mathrm{Li}_2(\mu^{\frac{M}{D}})).
\label{FNM}\end{eqnarray}
The transformation $(N,M) \rightarrow (\frac{N}{D},\frac{M}{D})$ does not change the first and the last summands. Hence this explicit
representation may also be used to derive the difference $F_{N,M}(t)-F_{\frac{N}{D},\frac{M}{D}}(t)$, that we have already computed
with different purely elementary strategy.

It remains to understand the case $M=0$. Since $F_{N,0}(t)=F_{1,0}(t)$, we need to consider
only the latter function. We have
\begin{align*}
\int_{0}^{1} \log |\theta(e^{2\pi i z}; e^{2\pi i\tau})|dx=\int_{0}^{1} \log |e^{-\frac{\pi i}{6}(\tau+3-\tilde \tau)+\pi i \tilde\tau z^2+\pi i(1-\tilde\tau)z}\theta(e^{2\pi i \tilde \tau z};
e^{2\pi i\tilde \tau })|dx,
\end{align*}
where $e^{2\pi i z}=te^{2\pi i x}$ and $\tilde \tau = -1/\tau$.
Here we have applied the modular transformation law for the theta function,
which we describe in terms of the uniform variables $\tau:=\omega_1/\omega_2, z:= u/\omega_2$:
\begin{equation}
\theta\left(e^{-2\pi{i}\frac{u}{\omega_1}};
e^{-2\pi{i}\frac{\omega_2}{\omega_1}}\right)
=e^{\pi{i}B_{2,2}(u;\mathbf{\omega})} \theta\left(e^{2\pi{i}\frac{u}{\omega_2}};
e^{2\pi{i}\frac{\omega_1}{\omega_2}}\right),
\label{modthetashort}\end{equation}
where $B_{2,2}$ is a multiple Bernoulli polynomial of the second order,
$$
B_{2,2}(u;\omega_1,\omega_2)=\frac{u^2}{\omega_1\omega_2}-\frac{u}{\omega_1}
-\frac{u}{\omega_2}+\frac{\omega_1}{6\omega_2}+\frac{\omega_2}{6\omega_1}+\frac{1}{2}.
$$

Since $e^{2\pi i \tilde \tau z}= t^{\tilde\tau}e^{2\pi i \tilde \tau x}$, from the previous computations
we have
$$
\int_{0}^{1} \log |\theta(t^{\tilde\tau}e^{2\pi i \tilde \tau x};e^{2\pi i\tilde \tau })|dx
=\frac{|\tau|^2}{4\pi \mathrm{Im}\,\tau}\left(\mathrm{Re}\, \left(\frac{\log t}{\tau}\right)\right)^2-\frac{\mathrm{Im}\,\tau}{\pi}\mathrm{Re}(\mathrm{Li}_2(\mu)),
$$
where by definition
$$
\mathrm{arg} \, \mu = \mathrm{Im}\,\left(-\frac{\log t}{\tau}\right) +\frac{\mathrm{Re}\, \left(-\frac{1}{\tau}\right)}{\mathrm{Im}\,\left(-\frac{1}{\tau}\right)}\log |e^{-\frac{\log t}{\tau}}|=\frac{\log|t|}{\mathrm{Im}\,\tau}.
$$
On the other hand,
$$
\int_{0}^{1} \log |e^{-\frac{\pi i}{6}(\tau+3-\tilde \tau)+\pi i \tilde\tau z^2+\pi i(1-\tilde\tau)z}|dx
=\frac{\pi}{6}\mathrm{Im}\,\tau +\frac{1}{2}\log|t|-\mathrm{Im}\,\left(\frac{\log^2 t}{4\pi\tau}\right).
$$
Combining two pieces, after some simplifications, we finally obtain the following result
\begin{equation}
F_{1,0}(t)=\frac{\log^2|t|}{4\pi \mathrm{Im}\,\tau}+\frac{1}{2}\log|t|+\frac{\pi}{6}\mathrm{Im}\,\tau-\frac{\mathrm{Im}\,\tau}{\pi}
\mathrm{Re}(\mathrm{Li}_2(e^{i\frac{\log|t|}{\mathrm{Im}\,\tau}})),
\label{F10}\end{equation}
which coincides with the formal substitution $(N,M)=(1,0)$ into the expression \eqref{FNM}
obtained earlier, since $\mu^{\frac{M}{D}}\big|_{M=0}=e^{i\frac{\mathrm{Re}\, \log \bar{t}}{\mathrm{Im}\,\tau}}$.
Thus, formula \eqref{FNM}
is well-defined and gives the correct answer for any $N,M \in \mathbb{Z}$.

We now return to our main problem of the integral evaluation in \eqref{sum_int}.
Since function \eqref{Hu} has equal number of theta functions in the numerator and
denominator, the second term in \eqref{FNM/D} cancels out automatically.
The balancing condition \eqref{balance} leads to the cancellation of the contributions from the
third term in \eqref{FNM/D}, since $\sum_{k=1}^r\log |t_k/w_k|=0.$
This means that the final expression does not depend on $D$,
\begin{equation}
\int_{0}^{1} \log |H(e^{2\pi i x (N+M\tau)})|d x =\int_{0}^{1} \log |H(e^{2\pi i x
(\frac{N}{D}+\frac{M}{D}\tau)})|d x.
\label{integralmodD}\end{equation}
Therefore, in this case we may assume that $D=(M,N)=1$. Then, we see that the balancing
condition removes contributions from the second and third terms in the evaluation \eqref{FNM}.
As a result, we come to the conclusion that the taken elliptic
hy\-per\-geo\-met\-ric series has the following radius of convergence
\begin{eqnarray}\nonumber &&  \makebox[-1em]{}
\log r_{c}^{-1}=\int_{0}^{1} \log |H(e^{2\pi i x (N+M\tau)})|d x
=\frac{1}{4\pi\mathrm{Im}\,\tau}\sum_{k=1}^r \left(\log^2|t_k|- \log^2|w_k|\right)+
\\ &&  \makebox[-2em]{}
+\frac{\mathrm{Im}\,\tau}{\pi|N+M\tau|^2} \mathrm{Re}\, \sum_{k=1}^r \left(\mathrm{Li}_2(e^{i\mathrm{Re}
\frac{(N+M\tau)\log \bar{w}_k}{\mathrm{Im}\,\tau}})-\mathrm{Li}_2(e^{i\mathrm{Re}\, \frac{(N+M\tau)\log \bar{t}_k}{\mathrm{Im}\,\tau}})
\right).
\label{r_fin}\end{eqnarray}
In view of the remark \eqref{bern} the second summand may be rewritten as a combination of
quadratic polynomials in fractional parts of $\arg z/2\pi$ with $z$ being the dilogarithm arguments.
As a result, we arrive to the following theorem.

\textbf{Theorem 5}. The radius of convergence of the elliptic hy\-per\-geo\-met\-ric series $_{r+1} E_{r}$ with $t_0=q=e^{2\pi i\chi (N+M\tau)}$, where $\chi $ is an arbitrary real irrational number, coprime integers
$N, M\in\ZZ$, $(N,M)=1$, and all of the $t_i, w_j$ not lying on the line $q^{\mathbb{R}}$
is given by the following explicit formula
\begin{eqnarray}\nonumber &&  \makebox[-1em]{}
\log r_{c}^{-1}=\sum_{k=1}^r \left(\frac{\log^2|t_k|- \log^2|w_k|}{4\pi\mathrm{Im}\,\tau}+
\frac{\pi\mathrm{Im}\,\tau}{|N+M\tau|^2} (\alpha_k(\alpha_k-1)-\beta_k(\beta_k-1))\right),
\\ &&  \makebox[1em]{}
\alpha_k:=\Biggl\{ \frac{\mathrm{Re}((N+M\tau)\log \bar{w}_k)}{2\pi\mathrm{Im}\,\tau}\Biggr\}, \quad
\beta_k:=\Biggl\{ \frac{\mathrm{Re}((N+M\tau)\log \bar{t}_k)}{2\pi\mathrm{Im}\,\tau}\Biggr\},
\label{r_fin_fin}\end{eqnarray}
where $\{x\}$ means the fractional part of $x$.\\

Consider now the radius of convergence of the well-poised elliptic hy\-per\-geo\-met\-ric series appearing
in the most interesting applications. Such series are defined by imposing the following set
of conditions
\begin{equation}
t_1w_1= \ldots =t_rw_r,
\label{WP}\end{equation}
so that there remains only $r$ free parameters, say, $t_1,\ldots, t_r$, and $w_1$ constrained by the
balancing condition $\prod_{j=2}^rt_j=\nu t_1^{\frac{r-2}{2}}w_1^{\frac{r}{2}},
\, \nu=\pm 1,$ in addition to
the qualitatively different variables $z, q$ and $p$.
One can check that all the terms of resulting $_{r+1}E_r$-series are elliptic functions of
parameters $t_j$ and $w_1$, i.e. they are invariant under the transformations
$t_j\to p^{n_j}t_j, w_1\to p^m w_1, \, n_j, m\in\mathbb{Z},$ respecting the balancing condition.
Namely,  one has the constraint $\sum_{j=2}^rn_j=\frac{r-2}{2}n_1 +\frac{r}{2}m$,
where for even $r$ and $\nu =1$ the discrete variables $n_1$ and $m$ can take
any integer value, whereas for even $r$ and $\nu =-1$ as well as for odd $r$ the sum
$n_1+m$ must be even. Since the invariance under the change $q\to pq$ is implemented
from the very beginning, only the variable $z$ does not lie on the elliptic curve with the
modulus $p$, i.e. it looks alien in the picture.

The constraints \eqref{WP} together with the
balancing condition force to vanish the first term in \eqref{r_fin_fin}.
In order to analyze the second term it is
convenient to use the following unique representation of the parameters
\begin{equation}
t_j=q^{h_j} e^{\varphi_j  \frac{2\pi \mathrm{Im}\,\tau}{N+M\bar\tau}}, \qquad
w_j=q^{\tilde{h}_j} e^{\tilde{\varphi}_j \frac{2\pi \mathrm{Im}\,\tau}{N+M\bar\tau}},
\quad h_j, \tilde{h}_j, \varphi_j, \tilde{\varphi}_j \in \mathbb{R}.
\label{hphi}\end{equation}
As we see, the variables $\varphi_j$ and $\tilde\varphi_j$ measure
how far the poles and zeros determined by the parameters $t_j$ and $w_j$
are from the line $q^\mathbb{R}$.
The variables $h_j, \tilde{h}_j, \varphi_j,$ and $\tilde{\varphi}_j$ are independent
(similar to the modulus and argument of a complex number) and satisfy the constraints
$$
\tilde{\varphi}_i+\varphi_i= \tilde{\varphi}_j+\varphi_j, \quad
\tilde{h}_i+h_i= \tilde{h}_j+h_j, \quad
\sum_{k=1}^r\tilde{\varphi}_k=\sum_{k=1}^r\varphi_k, \quad
\sum_{k=1}^r\tilde{h}_k =\sum_{k=1}^r h_k.
$$
The variables $h_j$ and $\tilde h_j$ do not give contribution to the real parts of interest
and the expression for the radius of convergence can be rewritten as follows
\begin{equation}
\log r_{c}^{-1}=\frac{\pi\mathrm{Im}\,\tau}{|N+M\tau|^2} \sum_{k=1}^r  (\{\tilde{\varphi}_k\}-\{\varphi_k\})(\{\tilde{\varphi}_k\}+\{\varphi_k\}-1).
\label{radius}\end{equation}
Evidently, if $\{\varphi_k\}=\varphi_k$, $\{\tilde{\varphi}_k\}=\tilde{\varphi}_k$ for all $k$,
then the sum on the right-hand side vanishes and $r_c=1$.

Consider a particular example of the parameter values when $\{\varphi_k\}\neq\varphi_k$ for some $k$.
Let us fix the following particular choice of the values of $\varphi_i$:
$$
\varphi_1=1+\frac{\varepsilon r}{2},\qquad \varphi_2=\ldots=\varphi_r =1-\varepsilon,
\quad \varepsilon >0.
$$
The indicated constraints on tilded variables can be resolved, which yields
$$
\tilde{\varphi}_j+\varphi_j=\frac{2}{r}\sum_{k=1}^r\varphi_k, \quad \tilde{\varphi}_1=1-\varepsilon \Big(1-\frac{2}{r}+\frac{r}{2}\Big),\quad \tilde{\varphi}_2=\ldots=\tilde{\varphi}_r=1+\frac{2\varepsilon}{r}.
$$
Consider now a particular form of the $\varepsilon$-variable
$$
\varepsilon = \frac{k+1}{\frac{r}{2} + \lambda},  \quad k \in \mathbb{Z}_{\ge 0},
$$
where $\lambda$ is a free real variable satisfying the constraints
$$
0< \lambda < 1-\frac{2}{r}, \quad r>2.
$$
For this particular choice we have the bounds
$$
\frac{k+1}{\frac{r}{2}+1-\frac{2}{r}} < \varepsilon<\frac{k+1}{\frac{r}{2}}.
$$
We require further $k+1 \leqslant \frac{r}{2}$, hence $\varepsilon < 1$ and $k+1 < \frac{r}{2}+\lambda$. Since $\lambda < 1$, we have $\frac{r}{2}-k\lambda > 0$. Adding to both sides $kr/2$ and dividing by $\lambda+ r/2$,
we come to the key constraint $\frac{\varepsilon r}{2} > k$, or $\varphi_1>k+1$.
The latter means that $\{\varphi_1\}=\frac{\varepsilon r}{2}-k$. Consequently, $\tilde{\varphi}_1 = 1-\varepsilon(\frac{r}{2}+1-\frac{2}{r})$. We also have $0<\varepsilon (-\lambda+(1-\frac{2}{r})) < 1$. Hence $\{\tilde{\varphi}_1\}=\tilde{\varphi}_1+k+1$. Also for $j>1$, $\{\varphi_j\}=\varphi_j$ and $\{\tilde{\varphi}_j\}=\tilde{\varphi}_j-1$. As a result, the straightforward computation yields
\begin{eqnarray} 
\log r_{c}^{-1}
=\frac{ \varepsilon\, \pi\mathrm{Im}\,\tau}{|N+M\tau|^2} (2\lambda + (\tfrac{2}{r}-1)(2k+4-r)).
\label{rc>1}\end{eqnarray}
The choice $k=[\frac{r-2}{2}]$ leads to the radius of convergence greater than $1$ for even $r$. Restricting
the value of $\lambda$ as $0< 2\lambda < 1-\frac{2}{r}$, we see that $r_c>1$ for odd $r$ as well.
Thus we have shown that there are highly nontrivial examples of the infinite well-poised elliptic
hy\-per\-geo\-met\-ric series with the radius of convergence bigger than 1. This is an
important fact, since it is the values $z=\pm 1$ that emerge most often in the applications (see below).
Note that $r_c>1$ for $q^b=p^M$ cases \eqref{roots} as well, provided $|R|<1$.

Suppose now that $r$ is even, $r=2m$ and $\nu=1$. Then we clear up the balancing condition by fixing
$t_r=t_1^{m-1}w_1^m/\prod_{j=2}^{r-1}t_j$ and consider the resulting convergent well-poised
series for fixed $z$ as analytical function of the free parameters:
$$
f(t_1,\ldots,t_{r-1},w_1)=\sum_{k=0}^\infty \prod_{n=0}^{k-1} h(n)\, z^k, \quad
h(n)=\prod_{j=1}^{r}
\frac{\theta(t_jq^n;p)}{\theta(\frac{t_1w_1}{t_j}q^n;p)}.
$$
This is a $p$-shift invariant function of all its variables
$$
f(t_1,\ldots, pt_j,\ldots,w_1)=f(t_1,\ldots, t_{r-1}, pw_1)=f(t_1,\ldots,t_{r-1},w_1),
$$
however, it is not an elliptic function, since it is not meromorphic.
E.g., as a function of $t_k$, for a fixed $k\in\{2,\ldots, r-1\},$ it has natural boundaries for
$t_k=t_1w_1q^\mathbb{R}p^n$ and $t_k=t_1^{m-2}w_1^{m-1}q^\mathbb{R}p^n/\prod_{j=2, \neq k}^{r-1}t_j$
with $n\in\mathbb{Z}$. As a function of $t_1$ or $w_1$, it has a different richer set of natural boundaries.

\section{The general case $t_0\neq q$}

Now we remove the constraint $t_0=q$ imposed earlier on the general $_{r+1}E_r$-series.
As in the previous situation, the balancing condition implies that $H(u)$ with the
argument $u=e^{2\pi i x}$ is an elliptic function of $x$ with the fundamental periods  $\{1,\tau\}$ whose
zeros located at the points $u=t_{i}^{-1}p^{\mathbb{Z}}$, $i=0,\ldots,r$, and poles --- at
the points $u=w_{j}^{-1}p^{\mathbb{Z}}$, $j=1,\ldots, r$. Consider any pair of integers $N$  and $M$,
such that the line from $x=0$ to $x=N+M\tau$ does not contain any zeros or poles of this
function, together with a real number $\chi  \in \mathbb{R}\setminus \mathbb{Q}$, and set
$$
q=e^{2\pi i \chi (N+M\tau)}.
$$
The reason why we restrict our attention to the irrational numbers $\chi $ only is precisely
the presence of the pole of $H(u)$ at the point $u=q^{-1}$. Since the singularity of $H(u)$
is of a particular nice type, the function $\log |H(e^{2\pi i x (N+M\tau)})\sin(\pi (x+\chi ))|$,
considered as a function on the unit interval, is continuous there. It follows that this function
fulfills the requirements of the Weyl equidistribution theorem, so that the next limit exists
for arbitrary irrational $\chi $:
\begin{eqnarray}\nonumber  &&
\lim_{n\to\infty} \frac{1}{n} \sum_{k=0}^{n-1} \log |H(e^{2\pi i k\chi  (N+M\tau)})\sin(\pi ((k+1)\chi ))|
\\ && \makebox[2em]{}
= \int_{0}^{1} \log |H(e^{2\pi i x (N+M\tau)})\sin(\pi (x+\chi ))| dx.
\label{sing}\end{eqnarray}

As follows from the considerations of Hardy and Littlewood \cite{HL}, this limit may be simplified
further for those irrational numbers, that satisfy the constraint on the denominators of their
Pad\'e approximants
\begin{equation}
\limsup_{k\to\infty} \frac{\log q_{k+1}}{q_k} = 0 \iff \liminf_{n \to \infty}
 |1-e^{2\pi i n\chi }|^{\frac{1}{n}}=1.
\label{HLcond}\end{equation}
And the latter condition holds for almost all irrational numbers. Therefore, for any such $\chi $
the $\sin$-function multiplier in \eqref{sing} can be removed, and we have
$$
r_{c}^{-1}:=\lim_{n \to \infty} \sqrt[n]{|H(1)\cdots H(q^{n-1})|}
= \exp\left(\int_{0}^{1} \log |H(e^{2\pi i x (N+M\tau)})|d x \right).
$$
Thus the computations of the previous section are applicable here as well --- it is just sufficient to extend
the summation
to the value of index $k=0$ and set $w_0=q$. We have the following result which extends that of Theorem 5
\begin{eqnarray}\nonumber &&  \makebox[-1em]{}
\log r_{c}^{-1}=\sum_{k=0}^r \left(\frac{\log^2|t_k|- \log^2|w_k|}{4\pi\mathrm{Im}\,\tau}+
\frac{\pi\mathrm{Im}\,\tau}{|N+M\tau|^2} (\alpha_k(\alpha_k-1)-\beta_k(\beta_k-1))\right),
\\ &&  \makebox[1em]{}
\alpha_k:=\Biggl\{ \frac{\mathrm{Re}((N+M\tau)\log \bar{w}_k)}{2\pi\mathrm{Im}\,\tau}\Biggr\}, \quad
\beta_k:=\Biggl\{ \frac{\mathrm{Re}((N+M\tau)\log \bar{t}_k)}{2\pi\mathrm{Im}\,\tau}\Biggr\},
\label{rt0_fin_fin}\end{eqnarray}
where $\{x\}$ means the fractional part of $x$. This formula in turn may be simplified a little more since $\mathrm{Re}((N+M\tau)\log \bar{q})=0$ implies $\alpha_0=0$.

Consider  how our analysis simplifies in the limit $p\to 0$. As $p$ approaches $0$, the expression
we used for $q=e^{2\pi i \chi (N+M\tau)}$ becomes meaningless for $N, M>0$.
Since we required $\frac{1}{2\pi i}\log q$ to lie on the line from $x=0$ to $x=N+M\tau$, the limit
of interest forces us to impose one of the following two conditions. Since Im$\,\tau\to\infty$, for  $M\neq 0$
our line becomes the vertical line and the sequence $q^n$ should ``walk'' along it.
Therefore Re$(\frac{1}{2\pi i}\log q)=0$, or $q=e^{-2\pi \chi }$ for some $\chi  \in \mathbb{R}^\times$.
For $M=0$ the sequence $q^n$ ``walks'' along the line Im$\, x=0$, i.e. $q=e^{2\pi i \eta}$ for some $\eta \in \mathbb{R}$. In the first case, we weaken the requirement
for $\frac{1}{2\pi i}\log t_i$  and $\frac{1}{2\pi i}\log w_i$, $i=1,\dots, r$, not to lie on the vertical ray, that contains $\frac{x}{2\pi i}\log q$ for $x\leqslant 0$. It is enough to consider the standard condition
\begin{equation}
t_j,\, w_j\neq  q^l, \quad l\in \mathbb{Z}_{\leqslant 0}, \quad j=0,\ldots, r.
\end{equation}
And in the second case we impose the condition for the zeros and poles not to lie on the unit interval.
The behaviour of the considered sequences (and hence of integrals) as $p$ approaches $0$ may then be described as follows. Since $H(u)$ itself has the limit
$$
H_0(u) := \lim_{p\to 0} H(u)=\frac{(1-ut_0)(1-ut_1)\dots(1-ut_r)}{(1-uq)(1-uw_1)\dots(1-uw_r)},
$$
in the case $M\neq 0$ the sequence $H_0(q^n)$ becomes bounded. The case $M=0$ and $t_0=q=e^{2\pi i \eta}$ for rational $\eta$ is also trivial. The case $M=0$ and $t_0=q=e^{2\pi i \eta}$ and irrational $\eta$ leads to the limit
$$
\lim_{n\to \infty} \frac{1}{n}\sum_{k=0}^{n-1} \log|H_0(q^k)| = \int_{0}^{1} \log|H_0(e^{2\pi i x})|dx.
$$

Finally, in the case $M=0$ and $t_0 \neq q = e^{2\pi i \eta}$, the variable $\eta$ is supposed to be irrational
by default due to the singularity of $H_0(u)$ at $u=q^{-1}$. For those irrational $\eta$, that satisfy the familiar constraint for the denominators of Pad\'e approximants
$$
\limsup_{k\to \infty} \frac{\log q_{k+1}}{q_k} = 0,
$$
the same limit exists
$$
\lim_{n\to\infty} \sqrt[n]{|H_0(1)\dots H_0(q^{n-1})|} = \exp\left(\int_{0}^{1} \log |H_0(e^{2\pi i x})|dx\right).
$$

Since we restrict our attention to the case $M=0$, the only $F_{N,M}(t)$ whose limiting behaviour is supposed to be examined is $F_{N,0}(t)=F_{1,0}(t)$. Set $C=[\frac{\log|t|}{2\pi \mathrm{Im}\,\tau}]$, where $[x]$ denotes the integer part of $x$. Using the remark \eqref{bern}
and the substitution $\{\frac{\log|t|}{2\pi \mathrm{Im}\,\tau}\}=\frac{\log|t|}{2\pi \mathrm{Im}\,\tau}-C$ in the equality
\eqref{F10}, we obtain an exact expression
$$
F_{1,0}(t) = (C+1)\log|t|-C(C+1)\pi\, \textrm{Im}\, \tau.
$$
For sufficiently large $\mathrm{Im}\, \tau$ the equality $C=0$ holds true, if $\log|t| > 0$, and $C=-1$, when $\log|t|<0$, and
we obtain
\begin{align*}
\lim_{p\to 0} F_{1,0}(t)
= \begin{cases}
      0, & \text{if}\ |t|\leq 1 \\
      \log|t|, & \text{if}\ |t|>1
    \end{cases}
\equiv\int_{0}^{1} \log |1-te^{2\pi ix}|dx.
\end{align*}
Therefore the appearance of the function $H_0(u)$ under the integral sign is fully justified.

Major applications of the elliptic hypergeometric series appear in the form of the following
very-well poised series \cite{spi:theta1}
\begin{eqnarray}
&&{}_{r+1}V_{r}(t_0;t_1,\ldots,t_{r-4};q,p):= \sum_{n=0}^\infty
\frac{\theta(t_0q^{2n};p)}{\theta(t_0;p)}\prod_{m=0}^{r-4}
\frac{\theta(t_m;p;q)_n}{\theta(qt_0t_m^{-1};p;q)_n}q^n,
\label{V-series}\end{eqnarray}
with the balancing condition
\begin{equation}
\prod_{k=1}^{r-4}t_k=\nu t_0^{(r-5)/2}q^{(r-7)/2}, \quad \nu=\pm 1,
\label{balancing}\end{equation}
where for odd $r$ it is assumed that $\nu=1$. This function can be obtained from
the general $_{r+1}E_r$ elliptic hypergeometric series
by imposing the well poisedness condition $qt_0=t_jw_j, \, j=1,\ldots, r$,
and additional five constraints
\begin{equation}
t_{r-3}=qt_0^{1/2}, \; t_{r-2}=-qt_0^{1/2}, \; t_{r-1}=q(t_0/p)^{1/2}, \; t_r=-q(pt_0)^{1/2},
\quad z=-1,
\label{vwpcond}\end{equation}
which result in the relation
$$
\prod_{j=r-3}^r \frac{\theta(t_jq^n;p)}{\theta(t_0q^{n+1}/t_j;p)}z^n
=\frac{\theta(t_0q^{2n};p)}{\theta(t_0;p)}q^n.
$$
Therefore the function $H(u)$ determining the ratio of sequential $_{r+1}E_r$-series coefficients,
$c_{n+1}/c_n=H(q^n)z$, acquires the form
\begin{equation}
H(u)=-q\frac{\theta(t_0q^2u^2;p)}{\theta(t_0u^2;p)}
\prod_{j=0}^{r-4} \frac{\theta(t_ju;p)}{\theta(qt_0u/t_j;p)}
\label{Hvwp}\end{equation}
with the balancing condition transformed to the form $q^4\prod_{k=0}^{r-4}t_k=\prod_{k=0}^{r-4}w_k$
where $w_0=q$ which reduces to relation \eqref{balancing} (remark: the choice $\nu=1$ for odd $r$ is a special convention \cite{spi:theta1}).

Since the taken restrictions on $t_j$ employ square roots of the variable $t_0$
and introduce additional dependence on the $q$-variable, it is not obvious that all terms in the
series \eqref{V-series} are elliptic functions of all its parameters (including $q$), i.e. that
they are invariant under the transformations $t_j\to q^{n_j}t_j,\, j=0,\ldots, r-4,$
and $q\to p^m q$ with the integers $n_j$ and $m$ respecting the balancing condition,
$\sum_{j=1}^{r-4}n_j=\frac{r-5}{2}n_0 +\frac{r-7}{2}m$.
However, it is the case, as can be checked directly.

In order to investigate convergence of infinite $_{r+1}V_r$-series we do not substitute constraints
\eqref{vwpcond} into the general formula for the radius of convergence \eqref{rt0_fin_fin},
but compute it anew using the function \eqref{Hvwp}. Since theta functions in the first ratio
contain $u^2$ in the arguments, we should compute $F_{N,M}(q^2t_0)-F_{N,M}(t_0)$
using the formula \eqref{FNM} with the choice $D=2$. However, we have
$\mathrm{Re}\big((N+M\tau)\log \overline{q^2t_0}\big)=\mathrm{Re}\big((N+M\tau)\log \overline{t_0}\big)$
and the $\textrm{Li}_2$-function parts from  \eqref{FNM} cancel out yielding
$$
F_{N,M}(q^2t_0)-F_{N,M}(t_0) = \frac{\log^2|q^2t_0|-\log^2|t_0|}{4\pi\mathrm{Im}\,\tau}-\frac{M-1}{2}\log\frac{|q^2t_0|}{|t_0|}.
$$
Contributions of other theta functions in \eqref{Hvwp} to the radius of convergence are
\begin{eqnarray} \nonumber &&
\sum_{k=0}^{r-4} (F_{N,M}(t_k)-F_{N,M}(w_k))
= \frac{1}{4\pi\mathrm{Im}\,\tau} \log\Big|\prod_{k=0}^{r-4}\frac{t_k}{w_k}\Big| \log|qt_0|
\\ && \makebox[-1em]{}
-\frac{M-1}{2}\log\big|\prod_{k=0}^{r-4}\frac{t_k}{w_k}\Big|
+\frac{\pi\mathrm{Im}\,\tau}{|N+M\tau|^2}\sum_{k=0}^{r-4} (\alpha_k(\alpha_k-1)-\beta_k(\beta_k-1)),
 \nonumber \end{eqnarray}
where $\alpha_k,\, \beta_k$ are defined in \eqref{rt0_fin_fin}. Adding these two expressions
and the contribution $\log|q|$ coming from the $q^n$-term in the $_{r+1}V_r$-series,
applying the balancing condition \eqref{balancing} and the parametrization \eqref{hphi},
we obtain
\begin{equation}
\log r_{c}^{-1}=M\log|q|+\frac{\pi\mathrm{Im}\,\tau}{|N+M\tau|^2} \sum_{k=0}^{r-4}  (\{\tilde{\varphi}_k\}-\{\varphi_k\})(\{\tilde{\varphi}_k\}+\{\varphi_k\}-1).
\label{rcvwp}\end{equation}
In \eqref{rcvwp} we have the constraints $\tilde \varphi_0=0$, following from
the convention $w_0=q$, and
$$
\varphi_0=\varphi_k+\tilde \varphi_k, \; k=1,\ldots, r-4,\quad
\sum_{k=0}^{r-4} \varphi_k=\sum_{k=0}^{r-4}\tilde \varphi_k.
$$
The latter equality follows from the relation $\prod_{k=0}^{r-4}t_k/w_k=q^{-4}$.

Assume now that $\{\varphi_k\}=\varphi_k, \, k=0,\ldots, r-4$, which is not forbidden by
the taken constraints. Then the sum of $\varphi_k$ and $\tilde \varphi_k$ depending
terms in \eqref{rcvwp} vanishes, and we obtain $\log r_{c}^{-1}=M\log|q|$. Because for $M>0$
we always have $|q|<1$, it follows that $r_c=|q|^{-M}>1$, i.e. the infinite $_{r+1}V_r$ very-well poised
elliptic hypergeometric series \eqref{V-series} does converge for the basic variable $q$ satisfying
the constraints \eqref{q_choice} and \eqref{HLcond} and the parameters $\varphi_k,\, \tilde \varphi_k$
(as fixed in  \eqref{hphi}) taken from the domain $0\leq \varphi_k,\, \tilde \varphi_k<1$,
$\varphi_k+\tilde \varphi_k<1,\, k=1,\ldots,r-4.$

The very-well poised part $q^n\theta(t_0q^{2n};p)/\theta(t_0;p)$ in the $_{r+1}V_r$-series
can be cancelled by a special choice of five parameters. Namely, one can give to $w_{r-7},\ldots, w_{r-4}$
the values $\pm qt_0^{1/2},$ $ q(t_0/p)^{1/2},$ $-q(pt_0)^{1/2}$
and choose $t_{r-8}=q$. After denoting $qt_0=t_1w_1$, there emerges exactly the series that we have considered
in the previous section with $z=-1$ and $r$ replaced by $r-9$.
Therefore the well-poised elliptic hypergeometric series with the radius
of convergence bigger than 1 \eqref{rc>1} yields another examples of converging infinite $_{r+1}V_r$-series.
Alternatively, one can take only one restriction $t_0=q$ (so that $t_kw_k=q^2$, $k=1,\dots,r-4$)
which makes the choice of $\chi$ as a rational number \eqref{roots} admissible, and then the
resulting $_{r+1}V_r$-series is summed as a geometric series for $|R|<1$.
All this opens up the question on whether one can represent the results of exact computation of
the elliptic beta integral and its multidimensional extensions \cite{spi:essays} as a
consequence of some identities for infinite elliptic hy\-per\-geo\-met\-ric series
(with appropriate restrictions on parameters) generalizing
the summation formulas for non-terminating very-well poised balanced $_8\varphi_7$-series
\cite{GR} and their multiple sum analogues.

Our analysis of the convergence of elliptic hy\-per\-geo\-met\-ric series used a special choice of the
parameter $q$ guaranteeing that the series coefficients do not have singularities beyond
the particular lines on the complex plane determined by the parameters $t_k$ and $w_k$.
For generic values of $p$ and $q$, for any choice of these parameters there emerge
series coefficients for which $q^n,\, n\to\infty,$ approach poles of the $H$-function arbitrarily close
and the question arises whether the two-dimensional version of the Weyl equidistribution theorem
is applicable in such a situation. This question remains open and requires a separate detailed investigation.

\section{Special cases}

We now restrict our attention to some special (the simplest) theta-hy\-per\-geo\-met\-ric series and
investigate, if there are conditions that imply their convergence. We shall take the only
constraint on the values of $q$ that $|q|\neq 1$. Considerations of the case $q=e^{2\pi i \chi }$,
$\chi \in\RR$, requires some deeper analysis in view of the recent results on the
bound $\lim\inf_{n\to\infty}|(q;q)_n|>0$ established in \cite{GKN} for $\chi $
given by some algebraic numbers, in particular, by the golden ratio.

$\bullet$   \, The $ _{0} E_0$-series:
$$
_{0} E_0\left( - ; q, p| z\right) = \sum_{n=0}^\infty \frac{z^n}{\theta(q ; p ; q)_n}, \quad 0 < |p| < 1.
$$
This series is interesting because it represents a direct theta-functional generalization
of one of the Euler $q$-exponential functions \cite{aar,GR}
$$
_{0} E_0\left(-;q, 0 | \, z\right) = \sum_{n=0}^\infty \frac{z^n}{(q; q)_n}=\frac{1}{(z;q)_\infty}, \quad |z|<1.
$$

The ratio of consecutive terms in $_{0} E_0$ is equal to $z/\theta(q^n; p)$ and we assume that
$q^k\neq p^l$ for $k, l \in \mathbb{Z}$, i.e. that the series is defined correctly. The standard theta function
transformation \eqref{trafo} implies that it is enough to check the case $0 < |q| <1$ to understand full picture.

Define a positive real number $\alpha = \log |q|/\log |p|$. For any positive integer $n$ define in addition the integer number $N_n=[n\alpha]$, so that the fractional part $\{n\alpha\}=n\alpha - N_n$. Then
$$
|\theta(q^n ; p)| = |\theta(p^{N_n} q^n p^{-N_n} ; p)|
= | q^{-nN_n}p^{\binom{N_n+1}{2}} \theta(q^n p^{-N_n} ; p)|.
$$
Now $|q^n p^{-N_n}|=|q|^n |p|^{-N_n}=e^{n\log |q|- [n\alpha] \log |p|}=|p|^{\{n\alpha\}}$. Also, since $|p|<1$, we have $\log |p| < 0$, hence $|p|< |p|^{\{n\alpha\}} \leqslant 1$. In addition, recall the trivial inequality $|1-u| \geqslant 1-|u|$. Hence we have the following sequence of bounds
\begin{eqnarray} \nonumber
|\theta(q^n ; p)| && \geqslant |q|^{-nN_n} |p|^{\binom{N_n+1}{2}}(1-|p|^{\{n\alpha\}})(1-|p|^{1-\{n\alpha\}})(|p|; |p|)_\infty^2\\ \nonumber
&& \geqslant |p|^{-N_n^2-N_n\{n\alpha\}+\binom{N_n+1}{2}}(1-|p|^{\{n\alpha\}})(1-|p|^{1-\{n\alpha\}})(|p|; |p|)_\infty^2\\ \nonumber
&& \geqslant |p|^{-N_n\{n\alpha\}-\binom{N_n}{2}}(1-|p|^{\{n\alpha\}})(1-|p|^{1-\{n\alpha\}})(|p|; |p|)_\infty^2
\\
&& \geqslant |p|^{-\binom{N_n}{2}}(1-|p|^{\{n\alpha\}})(1-|p|^{1-\{n\alpha\}})(|p|; |p|)_\infty^2.
\label{est1}\end{eqnarray}

Suppose now that $\alpha$ is a rational number, $\alpha=a/b$. Recall that $q^b\neq p^a$, but it is still possible that $|q|^b=|p|^a$ for some coprime integers $a$ and $b$. If $n$ is coprime to $b$, then $(b-1)/b \geqslant \{n\alpha\}=k/b \geqslant 1/b$. Hence $|p|^{1-1/b}\leqslant |p|^{\{n\alpha\}} \leqslant |p|^{1/b}$ and the bound \eqref{est1} may now be rewritten as
$$
|\theta(q^n ; p)| \geqslant |p|^{-\binom{N_n}{2}}(1-|p|^{1/b})(1-|p|^{1-(1-1/b)})(|p|; |p|)_\infty^2 > (1-|p|^{1/b})^2(|p|; |p|)_\infty^2.
$$

If $b|n$, i.e. $n=bk$, then $N_n=[bka/b]=ka$. Also $\{ka\}=0$, so that $|q^n p^{-N_n}|=|q^{kb} p^{-ka}|=1$.
Hence,
\begin{align*}
|\theta(q^{bk} ; p)| &= |p^{-\binom{ka}{2}} \theta(q^{kb} p^{-ka} ; p)|\\
&=|p^{-\binom{ka}{2}} 2\sin(\pi k \sigma)(pe^{2\pi i k\sigma}; p)_\infty (pe^{-2\pi i k\sigma}; p)_\infty|\\
&\geqslant|p^{-\binom{ka}{2}} 2\sin(\pi k \sigma)|(|p|; |p|)_\infty^2,
\end{align*}
where we have denoted $q^{b} p^{-a}=e^{2\pi i \sigma}$, $\sigma \in \mathbb{R} \setminus \mathbb{Q}$ (recall that if $\sigma \in \mathbb{Q}$, then the series is undefined). Now, the function $|\sin(\pi k \sigma)|$ is 1-periodic in $\sigma$ and, moreover, $|\sin(\pi k \sigma)|=|\sin(\pi \{k \sigma\})|=|\sin(\pi \{-k \sigma\})|$. Notice also, that $\{x\}+\{-x\}=1$ for $x \notin \mathbb{Z}$, hence $0 < \{x\} \leqslant \frac{1}{2} \iff \frac{1}{2} \leqslant \{-x\} < 1$. Choosing the appropriate sign and using the inequality $|\sin(\pi x)|\geqslant x/2$ for $0\leqslant x \leqslant \frac{1}{2}$, we finally obtain the bound
$$
|\theta(q^{bk} ; p)| \geqslant|p|^{-\binom{ka}{2}}\{\pm k\sigma\}(|p|; |p|)_\infty^2.
$$

Now consider only those irrational numbers, that have the property
$$
\exists C(\sigma)>0: \inf_{k\in \mathbb{Z}_{>0}} |p|^{-\binom{ka}{2}}\{\pm k\sigma\} > C(\sigma).
$$
For example, this can be any irrational algebraic number, or $\pi$,
 or any other $\sigma$, such that the denominators of its Pad\'e approximants satisfy the inequality
$$
\limsup_{k\to \infty} \frac{\log q_{k+1}}{q_k^2} < -\frac{a^2}{2} \log |p|.
$$
Then the series converges for
$$
|z|<\min(C(\sigma),  (1-|p|^{1/b})^2)(|p|; |p|)_\infty^2.
$$

Suppose now that $\alpha$ is irrational. Then
\begin{align*}
|\theta(q^n ; p)| &\geqslant |p|^{-\binom{N_n}{2}}(1-|p|^{\{n\alpha\}})(1-|p|^{\{-n\alpha\}})(|p|; |p|)_\infty^2
\\  &
\geqslant |p|^{-\binom{N_n}{2}}(-\{n\alpha\}\log|p||p|^{\{n\alpha\}})(-\{-n\alpha\}\log|p||p|^{\{-n\alpha\}})(|p|; |p|)_\infty^2\\
&>|p|^{-\binom{n\alpha-1}{2}}\{n\alpha\}\{-n\alpha\}|p|\log^2|p|(|p|; |p|)_\infty^2,
\end{align*}
where in the last line we used a trivial inequality $N_n > n\alpha-1$. Consider again only those irrational numbers, that have the property
$$
\exists C(\alpha)>0: \inf_{n\in \mathbb{Z}_{>0}} |p|^{-\binom{n\alpha-1}{2}}\{n\alpha\}\{-n\alpha\} > C(\alpha),
$$
for example, any irrational algebraic number or $\pi$.
Then such $\alpha$ fits and the series converges in the domain
$$
|z|<C(\alpha)|p|\log^2|p|(|p|; |p|)_\infty^2.
$$

$\bullet$   \, The $ _{1} E_0$-series:
$$
_{1} E_0\left(\begin{matrix} t_0 \\ - \end{matrix}; q, p \, \bigg| \, z\right) = \sum_{n=0}^\infty \frac{\theta(t_0 ; p ; q)_n}{\theta(q ; p ; q)_n}z^n, \quad 0 < |p| < 1.
$$
This is a theta-functional analogue of the series standing on the left-hand side of
the $q$-binomial theorem \cite{aar,GR}
$$
_{1} E_0\left(\begin{matrix} t_0 \\ - \end{matrix}; q, 0 \, \bigg| \, z\right) =
\sum_{n=0}^\infty \frac{(t_0;q)_n}{(q;q)_n}z^n=\frac{(t_0 z;q)_\infty}{(z;q)_\infty}, \quad 0<|q|<1, \; |z|<1.
$$
The ratio of consecutive terms in $_{1} E_0$ is equal to $z\theta(t_0q^{n-1};p)/\theta(q^n;p)$
which is well defined for  $q^k\neq p^l$ for $k,l \in \mathbb{Z}$.   Again, set  $\alpha=\log|q|/\log|p|$. For any positive integer $n$ define $N_n=[n\alpha]$, so that $|q^n p^{-N_n}|=|p|^{\{n\alpha\}}$. Then from
  \eqref{trafo} we have the following equality
$$
\left|\frac{\theta(t_0q^{n-1};p)}{\theta(q^n;p)}\right|=|qt_0^{-1}|^{N_n} \left|\frac{\theta(t_0q^{n-1}p^{-N_n};p)}{\theta(q^np^{-N_n};p)}\right|
$$
We proceed as in the previous example . Suppose that $\alpha$ is a rational number. Recall that $q^b \neq p^a$,
but it is possible that $|q|^b=|p|^a$ for some coprime integers $a,b$. If $n$ is coprime to $b$, then
$(b-1)/b \geqslant \{n\alpha\} \geqslant 1/b$. Hence
$|p|^{1-1/b} \leqslant |p|^{\{n\alpha\}} \leqslant |p|^{1/b}$. In addition, the absolute value of theta function
is clearly bounded from above in the fundamental domain by some constant $S$. Thus we come to the following bound
$$
\left|\frac{\theta(t_0q^{n-1}p^{-N_n};p)}{\theta(q^np^{-N_n};p)}\right| \leqslant
 \frac{|\theta(t_0q^{n-1}p^{-N_n};p)|}{(1-|p|^{1/b})^2 (|p|;|p|)_\infty^2}
\leqslant \frac{ S}{(1-|p|^{1/b})^2 (|p|;|p|)_\infty^2}.
$$

If $b|n$, i.e. $n=bk$, then $N_n=ka$. Also $\{k\alpha\}=0$, so that $|q^n p^{-N_n}|=|q^{kb}p^{-ka}|=1$.
Hence, we have
$$
\left|\frac{\theta(t_0q^{kb-1};p)}{\theta(q^{kb};p)}\right| \leqslant \frac{1}{2}|qt_0^{-1}|^{ka}\frac{|\theta(t_0q^{-1}e^{2\pi i \sigma k};p)|}{|\sin(\pi k\sigma)| (|p|;|p|)_\infty^2}.
$$
Here we write $q^bp^{-a}=e^{2\pi i \sigma}$, for $\sigma \in \mathbb{R}\setminus \mathbb{Q}$. We now choose such $t_0$ that
$$
\exists C(\sigma, t_0) >0: \inf_{k\in\mathbb{Z}_{>0}} |qt_0^{-1}|^{-ka} |\{\pm k\sigma\}| > C(\sigma, t_0)
$$
and such $t_0$ exists among ``normal'' irrational numbers. As an example, we may choose $t_0 = eq$ or $2q$, with an appropriate irrational number $\sigma$. Then any such pair ($t_0, \sigma$) fits, and we have the bound
$$
\left|\frac{\theta(t_0q^{kb-1};p)}{\theta(q^{kb};p)}\right| \leqslant \frac{|\theta(t_0q^{-1}e^{2\pi i \sigma k};p)|}{C(\sigma, t_0)(|p|;|p|)_\infty^2} \leqslant \frac{S}{C(\sigma, t_0)(|p|;|p|)_\infty^2}.
$$
In other words, the series converges for
$$
|z| < \min ( (1-|p|^{1/b})^2, C(\sigma, t_0))S^{-1}(|p|;|p|)_\infty^2.
$$

Suppose now that $\alpha$ is irrational. We can write $q^n p^{-N_n}=|p|^{\{n\alpha\}} e^{2\pi i r_n}$
for some real numbers $r_n$. The absolute value of theta function in the numerator is bounded from above
and, also, $1-|p|^{\{n\alpha\}} \geqslant |p|^{\{n\alpha\}}(-\log|p|)\{n\alpha\}$. Hence
$$
\left|\frac{\theta(t_0q^{n-1};p)}{\theta(q^n;p)}\right| \leqslant  \frac{|qt_0^{-1}|^{N_n}|\theta(t_0q^{-1} |p|^{n\alpha} e^{2\pi i r_n};p)|}{\{n\alpha\}\{-n\alpha\}|p|\log^2|p| (|p|;|p|)_\infty^2} \leqslant  \frac{|qt_0^{-1}|^{N_n} S}{\{n\alpha\}\{-n\alpha\}|p|\log^2|p| (|p|;|p|)_\infty^2}.
$$
We now choose such $t_0$ that
$$
\exists C(\alpha, t_0) >0: \inf_{n\in\mathbb{Z}_{>0}} |qt_0^{-1}|^{-N_n} |\{n\alpha\}\{-n\alpha\}| > C(\alpha, t_0).
$$
Then any such pair ($t_0, \alpha$) fits, and the series converges for
$$
|z| < C(\alpha, t_0) S^{-1}|p|\log^2|p|(|p|;|p|)_\infty^2.
$$

\medskip

{\bf Acknowledgments.} This study has been partially funded within the framework
of the HSE University Basic Research Program.

\bibliographystyle{unsrt}

\end{document}